\documentclass[reqno, 11pt,fleqn]{amsart} 
\usepackage{graphicx}
\usepackage{amsmath,amssymb,amsthm, mathrsfs, mathtools,bbm}
\usepackage[usenames,dvipsnames]{xcolor}

\usepackage{enumerate}
\usepackage{tabularx}
\usepackage{ltablex}
\usepackage{booktabs}
\usepackage{colortbl}
\usepackage{bm}
\usepackage{esvect}
\usepackage{hyperref}
\hypersetup{
    colorlinks,
    linkcolor={red!60!black},
    citecolor={green!60!black},
    urlcolor={blue!60!black}
}
\usepackage[T1]{fontenc}
\usepackage{lmodern}
\usepackage[babel]{microtype}
\usepackage[english]{babel}
\usepackage{comment}

\usepackage{enumitem}
\usepackage{arydshln}

\usepackage[normalem]{ulem} 
\usepackage{dsfont}

\usepackage{geometry}
\graphicspath{{images/}} 
\usepackage{theoremref}

\usepackage{float}

\theoremstyle{plain}
\newtheorem{thm}{Theorem}[section]
\newtheorem{theorem}[thm]{Theorem}
\newtheorem*{theorem*}{Theorem}
\newtheorem{prop}[thm]{Proposition}
\newtheorem{obs}[thm]{Observation}
\newtheorem{claim}[thm]{Claim}
\newtheorem{cor}[thm]{Corollary}
\newtheorem{conj}[thm]{Conjecture}
\newtheorem{lemma}[thm]{Lemma}

\theoremstyle{definition}
\newtheorem{definition}[thm]{Definition}
\newtheorem{rem}[thm]{Remark}

\numberwithin{equation}{section}
\numberwithin{figure}{section}

\usepackage[hypcap=false]{caption}
\DeclarePairedDelimiter{\parens}{(}{)}
\DeclarePairedDelimiter{\set}{\{}{\}}

\DeclarePairedDelimiter\size{\lvert}{\rvert}   

\DeclareMathOperator{\claw}{Claw}

\renewcommand{\leq}{\leqslant}
\renewcommand{\geq}{\geqslant}

\renewcommand{\phi}{\varphi}
 
\newcommand{\sym}{{\text{sym}}}

\newcommand{\cM}{\ensuremath{\mathcal{M}}}

\newcommand{\cS}{\ensuremath{\mathcal{S}}} 
\newcommand{\cT}{\ensuremath{\mathcal{T}}}

\newcommand{\cV}{\ensuremath{\mathcal{V}}}
\newcommand{\cH}{\ensuremath{\mathcal{H}}}

\newcommand{\nto}{\nrightarrow}

\newcommand{\Burr}{Burr}
\newcommand{\Chvatal}{Chv{\'a}tal}

\newcommand{\Dudek}{Dudek}
\newcommand{\Erdos}{Erd\H{o}s}
\newcommand{\Folkman}{Folkman}
\newcommand{\Grinshpun}{Grinshpun}

\newcommand{\Han}{H{\`a}n}
\newcommand{\Kohayakawa}{Kohayakawa}
\newcommand{\Kreuter}{Kreuter}

\newcommand{\Lovasz}{Lov{\'a}sz}
\newcommand{\Luczak}{{\L}uczak}
\newcommand{\Nesetril}{Ne\v{s}et\v{r}il}

\newcommand{\Rodl}{R\"{o}dl}

\newcommand{\Rucinski}{Ruci\'{n}ski}
\newcommand{\Siggers}{Siggers}
\newcommand{\Szabo}{Szab\'o}
\newcommand{\Szekeres}{Szekeres}

\newcommand{\Turan}{Tur\'{a}n}
\newcommand{\Zumstein}{Zumstein}
\newcommand{\Zurcher}{Z\"urcher}

\newcommand{\FGLPS}{Fox, \Grinshpun{}, Liebenau, Person, and \Szabo{}}
\newcommand{\BEL}{\Burr{}, \Erdos{}, and \Lovasz{}}

\renewcommand{\mathbf}[1]{\bm{\mathcal{#1}}}

\setlist[enumerate,1]{label={(\roman*)}}
\setlist[enumerate,2]{label={(\alph*)}}

\title[Senders for asymmetric tuples of cliques in Ramsey theory]{On the use of senders \\ for asymmetric tuples of cliques in Ramsey theory}

\author[S.Boyadzhiyska]{Simona Boyadzhiyska}
\address[SB]{School of Mathematics, University of Birmingham, Edgbaston, Birmingham, B15 2TT, UK. Much of this research was done when the author was affiliated with the Institut f\"ur Mathematik, Freie Universit\"at Berlin, Berlin, Germany.}
\email{s.s.boyadzhiyska@bham.ac.uk}

\author[T.Lesgourgues]{
 Thomas Lesgourgues}  
\address[TL]{Combinatorics and Optimization Department, University of Waterloo, Waterloo, Ontario N2L 3G1, Canada. Much of this research was done when the author was affiliated with the School of Mathematics and Statistics, University of New South Wales, Sydney, Australia.}
\email{tlesgourgues@uwaterloo.ca}

\usepackage{tikz}
\usepackage{calc}
\usetikzlibrary{calc,angles,quotes,decorations.pathreplacing,calligraphy,backgrounds}

\tikzstyle{vertex}=[circle, draw, fill=black, inner sep=0pt, minimum size=6pt]
\tikzstyle{smallvertex}=[circle, draw, fill=black, inner sep=0pt, minimum size=4pt]

\newcommand{\vertex}{\node[vertex]}

\newcommand{\widthedge}{1}

\tikzset
{%
  pics/RSgraph/.style={%
    code={%
    
    \vertex (a) at (-3,0) [label=]{};
    \vertex (b) at (-1,0) [label=]{};
    \vertex (c) at (1,0) [label=]{};
    \vertex (d) at (3,0) [label=]{};

    \draw (0,0) ellipse (4 and 1);
    
    \draw (2.5,-4) arc
        [
            start angle=0,
            end angle=-180,
            x radius = 2.5,
            y radius = 1
        ] ;  
        
    \path
    (-4,0) edge (-2.5,-4)
    (4,0) edge (2.5,-4);
    
    \path
    (a) edge (-2,-2)
    (a) edge (-1.1,-2.4)
    (b) edge (-0.8,-2)
    (c) edge (0,-2)
    (c) edge (0.5,-2)
    (c) edge (1,-2)
    (d) edge (1.8,-3);
    }},
}

\newcommand{\tikzRS}[1]{
    \begin{tikzpicture}[x=#1 cm, y=#1 cm]

    \vertex (x) at (0,5) [label=20:$x$]{};
    \pic (test) at (0,0) {RSgraph} ;
    \path 
    (x) edge (testa)
    (x) edge (testb)
    (x) edge (testc)
    (x) edge (testd);

    \draw (140:4.5) node{$F$};
    
    \end{tikzpicture}
}

\newcommand{\tikzRSFull}[1]{
    \begin{tikzpicture}[x=#1 cm, y=#1 cm]

    \vertex (x) at (6,6) [label=20:$x_0$]{};
    \foreach \i in {1, 2, 3} {
        \pgfmathsetmacro{\sh}{6*(\i-1)}
        \begin{scope}[shift={(\sh,0)}]
          \pic (test\i) at (0,0) {RSgraph} ;
        \end{scope}   
        \path 
        (x) edge (test\i a)
        (x) edge (test\i b)
        (x) edge (test\i c)
        (x) edge (test\i d);        
            }
    \draw (140:4.5) node{$S$};
    
    \end{tikzpicture}
}

\newcommand{\tikzRSHyper}[1]{
    \begin{tikzpicture}[x=#1 cm, y=#1 cm]

    \foreach \i in {1, 2, ..., 10}{
        \pgfmathsetmacro{\j}{2*\i}
        \vertex (v\i) at (\j,0) [label=]{};
        }
    \draw (5,0) ellipse (4 and 1);
    \draw (11,0) ellipse (4 and 1);
    \draw (17,0) ellipse (4 and 1);
    
    \draw (90:1) node{$\mathbf{H}$};
    
    \end{tikzpicture}
}

\newcommand{\tikzjoinS}[1]{
\begin{tikzpicture}[x=#1 cm, y=#1 cm]

    \clip (-1.4,-1.2) rectangle (3.5,1.2);    
    \draw (180:0.75) node{$G_1$};
    
    \vertex (a_1) at (315:1) [label=-90:$a_1$]{};
    \vertex (b_1) at (0,0) [label=180:$b_1$]{};
    \vertex (c_1) at (45:1) [label=90:$c_1$]{};
    
    \path 
    (a_1) edge[line width= \widthedge] (b_1) 
    (b_1) edge[line width= \widthedge] (c_1);

    \begin{scope}[on background layer]
        \fill[gray!50] (a_1.center) to [out=210, in=150, looseness=5.4](c_1.center) to (b_1.center);
    \end{scope}
    
    \begin{scope}[shift={(1,0)}]

        \draw (0:1) node{$G_2$};
        
        \vertex (a_2) at (315:1) [label=-90:$a_2$]{};
        \vertex (b_2) at (0,0) [label=0:$b_2$]{};
        \vertex (c_2) at (45:1) [label=90:$c_2$]{};
        
        \path 
        (a_2) edge[line width= \widthedge] (b_2) 
        (b_2) edge[line width= \widthedge] (c_2);
        
        \begin{scope}[on background layer]
            \fill[gray!50] (a_2.center) to [out=-30, in=30, looseness=4.3](c_2.center) to (b_2.center);
        \end{scope}        
    \end{scope}
    
    \path[-latex] ([shift={(0.85,0)}]315:1) edge[line width= \widthedge, dashed] ([shift={(0.1,0)}]a_1);
    \path[-latex] ([shift={(0.85,0)}]0,0) edge[line width= \widthedge, dashed] ([shift={(0.1,0)}]b_1);
    \path[-latex] ([shift={(0.85,0)}]45:1) edge[line width= \widthedge, dashed] ([shift={(0.1,0)}]c_1);
    
\end{tikzpicture}}

\newcommand{\tikzAttachDet}[1]{
\begin{tikzpicture}[x=#1 cm, y=#1 cm]

    \clip (-1.4,-0.25) rectangle (2.5,1.25);   
    
    \draw (-0.75,0.5) node{$G$};
    \vertex (a) at (0,1) [label=]{};
    \vertex (b) at (0,0) [label=]{};
    \path 
    (a) edge[line width= \widthedge] node[left]{$f$}(b); 

    \draw (1.75,0.5) node{$R$};
    \vertex (c) at (1,1) [label=]{};
    \vertex (d) at (1,0) [label=]{};
    \path 
    (c) edge[line width= \widthedge] node[right]{$e$} (d);

    \path[-latex] ([shift={(0.85,0)}]c) edge[line width= \widthedge, dashed] ([shift={(0.1,0)}]a);
    \path[-latex] ([shift={(0.85,0)}]d) edge[line width= \widthedge, dashed] ([shift={(0.1,0)}]b);

    \begin{scope}[on background layer]
        \fill[gray!50] (a.center) to [out=150, in=210, looseness=5](b.center);
        \fill[gray!50] (c.center) to [out=30, in=-30, looseness=5](d.center);        
    \end{scope}

\end{tikzpicture}}

\newcommand{\tikzJoinSignal}[1]{
\begin{tikzpicture}[x=#1 cm, y=#1 cm]

    \clip (-2.5,-3.5) rectangle (4,1.5);    
    
    \draw (-1,-1) node{$G$};
    
    \vertex (a) at (1,1) [label=]{};
    \vertex (b) at (0,0) [label=]{};
    \vertex (c) at (0,-2) [label=]{};
    \vertex (d) at (1,-3) [label=]{};    
    
    \path 
    (a) edge[line width= \widthedge] node[left]{$a$} (b) 
    (c) edge[line width= \widthedge] node[left,shift=({-0.1,0})]{$b$} (d);
    
    \draw (3,-1) node{$S$};
    
    \vertex (u) at (3,1) [label=]{};
    \vertex (v) at (2,0) [label=]{};
    \vertex (w) at (2,-2) [label=]{};
    \vertex (x) at (3,-3) [label=]{};    
    
    \path 
    (u) edge[line width= \widthedge] node[right]{$e$} (v) 
    (w) edge[line width= \widthedge] node[right,shift=({0.1,0})]{$f$} (x);
    
    \path[-latex] ([shift={(-0.1,0)}]u) edge[line width= \widthedge, dashed] ([shift={(0.1,0)}]a);
    \path[-latex] ([shift={(-0.1,0)}]v) edge[line width= \widthedge, dashed] ([shift={(0.1,0)}]b);    
    \path[-latex] ([shift={(-0.1,0)}]w) edge[line width= \widthedge, dashed] ([shift={(0.1,0)}]c);
    \path[-latex] ([shift={(-0.1,0)}]x) edge[line width= \widthedge, dashed] ([shift={(0.1,0)}]d);

    \begin{scope}[on background layer]
        \fill[gray!50] (a.center) to [out=150, in=210, looseness=3.2](d.center) to (c.center) to [out=135, in=225, looseness=0.6](b.center);
        \fill[gray!50] (u.center) to [out=-45, in=45, looseness=1.1](x.center) to (w.center) to [out=45, in=-45, looseness=0.6](v.center);        
    \end{scope}
    
\end{tikzpicture}}

\newcommand{\tikzSTwoNew}[1]{
\begin{tikzpicture}[x=#1 cm, y=#1 cm]
    \vertex (x1) at (0,0) [label=90:$v$]{};
    \vertex (u1) at (-1,-2) [label=-90:$u$]{};
    \vertex (v1) at (1,-2) [label=-90:$x$]{};
    \vertex (v2) at (2,0) [label=90:$y$]{};
    
    \path
    (u1) edge[line width= 1.5]  node[right]{$e_1$}  (x1)
    (x1) edge[line width= \widthedge]  node[left]{$f_1$} node[right]{$e_2$}(v1)
    (v1) edge[line width= 1.5]  node[left]{$f_2$} node[right]{$h$}(v2)
    
    (u1) edge[line width= \widthedge, out = -20, in = -160, looseness = 1,dashed]  (v1)
    (x1) edge[line width= \widthedge, out = 20, in = 160, looseness = 1,dashed]   (v2);

    \begin{scope}[on background layer]
        \fill[gray!10] (u1.center) to [out = -20, in = -160, looseness = 1](v1.center) to (x1.center);
        \fill[gray!30] (x1.center) to [out = 20, in = 160, looseness = 1](v2.center) to (v1.center); 
    \end{scope}
    
    \draw (0,-1.5) node{\Large $\bm{S_1}$};
    \draw (1,-0.5) node{\Large $\bm{S_2}$};
    
\end{tikzpicture}}

\newcommand{\tikzSymmetricSVariant}[1]{
\begin{tikzpicture}[x=#1 cm, y=#1 cm]
        
    \vertex (a) at (0,1) [label=90:$a$]{};
    \vertex (b) at (0,0) [label=0:$b$]{};
    \vertex (c) at (0,-1) [label=-90:$c$]{};

    \path
    (a) edge[line width= \widthedge]  node[left]{$e_1$}  node[right]{$f_2$}(b)
    (b) edge[line width= \widthedge]  node[left]{$f_1$} node[right]{$e_2$}(c)
    (a.center) edge[line width= \widthedge, out = 180, in = 180, looseness = 1,dashed]  (c.center)
    (a.center) edge[line width= \widthedge, out = 0, in = 0, looseness = 1,dashed]  (c.center);

    \begin{scope}[on background layer]
        \fill[gray!25] (a.center) to [out = 180, in = 180, looseness = 1](c.center) to (b.center);
        \fill[gray!50] (a.center) to [out = 0, in = 0, looseness = 1](c.center) to (b.center);
    \end{scope}
    
    \draw (-0.4,0) node{\Large $\bm{S_1}$};
    \draw (0.4,0) node{\Large $\bm{S_2}$};

\end{tikzpicture}}

\newcommand{\tikzclawComplement}[1]{
\begin{tikzpicture}[x=#1 cm, y=#1 cm]
    \begin{scope}[on background layer]
        \draw [black, fill=black!30] (0,0) ellipse (1 and 4);
    \end{scope}
    \draw (0.5,4.5) node{\Large $H\cong K_7$};
    \foreach \i in {-3,-2,...,3} {
        \vertex (u\i) at (0,\i) [label=]{};
    }
    \vertex (z) at (0,0) [label=0:\Large $v_4$]{};
    \vertex (v1) at (0,3) [label=0:\Large $v_1$]{};
    \vertex (v2) at (0,2) [label=0:\Large $v_2$]{};
    \vertex (v3) at (0,1) [label=0:\Large $v_3$]{};
    
    \vertex (x) at (-4,0) [label=90:\Large $x$]{};
    \vertex (y) at (-6,0) [label=90:\Large $y$]{};
    
    \begin{scope}[on background layer]
        \fill[red!10] (y.center) to [out = 60, in = 180, looseness = 1](u3.center) to (x.center);
        \fill[red!20] (y.center) to [out = 60, in = 180, looseness = 1](u2.center) to (x.center);
        \fill[red!30] (y.center) to [out = 60, in = 180, looseness = 1](u1.center) to (x.center);
    \end{scope}

    \path (x) edge[line width= \widthedge]  (y) ; 
   
    \path 
    (x) edge[line width= \widthedge, red] (u1)
    (x) edge[line width= \widthedge, red] (u2)
    (x) edge[line width= \widthedge, red] (u3);
    
    \path
     (y) edge[line width= 1, red, out = 60, in = 180, looseness = 1] (u3)
     (y) edge[line width= 1, red, out = 60, in = 180, looseness = 1] (u2)
     (y) edge[line width= 1, red, out = 60, in = 180, looseness = 1] (u1);
    
    \path 
    (x) edge[line width= \widthedge,  blue] (z);    
    
    \draw[red] (-2,2.5) node{$S$};
    
    \draw [ultra thick, decorate,
        decoration = {calligraphic brace,amplitude=5pt}] (2,3) --  (2,1) node[pos=0.5,right = 5pt]{\Large $d=3$};
    
\end{tikzpicture}}

\newcommand{\tikzclawComplementInduction}[1]{
\begin{tikzpicture}[x=#1 cm, y=#1 cm]

    \begin{scope}[on background layer]
        \draw [black, fill=black!30] (0,0) ellipse (1 and 4);
    \end{scope}

    \foreach \i in {-3,-2,...,3} {
        \vertex (u\i) at (0,\i) [label=]{};
    }
    \vertex (z) at (0,1) [label=0:\Large $v_3$]{};
    
    \vertex (x) at (-4,0) [label=90:\Large $x$]{};
    \vertex (y) at (-6,0) [label=90:\Large $y$]{};
    
    \vertex (v1) at (0,3) [label=0:\Large $v_1$]{};
    \vertex (v2) at (0,2) [label=0:\Large $v_2$]{};
    
    \begin{scope}[on background layer]
        \fill[black!10] (y.center) to [out = 60, in = 180, looseness = 1](u3.center) to (x.center);
        \fill[black!15] (y.center) to [out = 60, in = 180, looseness = 1](u2.center) to (x.center);
    \end{scope}

    \path (x) edge[line width= \widthedge]  (y) ; 
   
    \path 
    (x) edge[line width= \widthedge] (u2)
    (x) edge[line width= \widthedge] (u3);
    
    \path
     (y) edge[line width= 1, out = 60, in = 180, looseness = 1] (u3)
     (y) edge[line width= 1, out = 60, in = 180, looseness = 1] (u2);
    
    \path 
    (x) edge[line width= \widthedge] (z); 

    \draw (2,0) node{\Huge $\Longrightarrow$};

    \begin{scope}[shift={(9,0)}]
        
        \begin{scope}[on background layer]
            \draw [black, fill=black!30] (0,0) ellipse (1 and 4);
        \end{scope}

        \foreach \i in {-3,-2,...,3} {
            \vertex (u\i) at (0,\i) [label=]{};
        }
        \vertex (z) at (0,0) [label=0:\Large $v_4$]{};
        
        \vertex (x) at (-4,0) [label=90:\Large $x$]{};
        \vertex (y) at (-6,0) [label=90:\Large $y$]{};
        \vertex (v1) at (0,3) [label=0:\Large $v_1$]{};
        \vertex (v2) at (0,2) [label=0:\Large $v_2$]{};    
        \vertex (v3) at (0,1) [label=0:\Large $v_3$]{};    
        
        \begin{scope}[on background layer]
            \fill[black!10] (y.center) to [out = 60, in = 180, looseness = 1](u3.center) to (x.center);
            \fill[black!15] (y.center) to [out = 60, in = 180, looseness = 1](u2.center) to (x.center);
            \fill[red!30] (y.center) to [out = 60, in = 180, looseness = 1](u1.center) to (x.center);
        \end{scope}
    
        \path (x) edge[line width= \widthedge]  (y) ; 
       
        \path 
        (x) edge[line width= \widthedge] (u1)
        (x) edge[line width= \widthedge] (u2)
        (x) edge[line width= \widthedge] (u3);
        
        \path
         (y) edge[line width= 1, out = 60, in = 180, looseness = 1] (u3)
         (y) edge[line width= 1, out = 60, in = 180, looseness = 1] (u2)
         (y) edge[line width= 1, red, out = 60, in = 180, looseness = 1] (u1);
        
        \path 
        (x) edge[line width= \widthedge, blue] (z);    
        
    \end{scope}

\end{tikzpicture}}

\newcommand{\tikzSScomplement}[1]{
\begin{tikzpicture}[x=#1 cm, y=#1 cm]
    \vertex (x1) at (0,0) [label=90:$b$]{};
    \vertex (u1) at (-1,-2) [label=-90:$a$]{};
    \vertex (v1) at (1,-2) [label=-90:$c$]{};
    \vertex (v2) at (2,0) [label=90:$d$]{};
    
    \path
    (u1) edge[line width= 1.5]   node[right]{$e$} (x1)
    (x1) edge[line width= \widthedge]  node[left]{$f$} node[right]{$e^c$}(v1)
    (v1) edge[line width= 1.5]  node[left]{$f^c$}(v2)
    
    (u1) edge[line width= \widthedge, out = -20, in = -160, looseness = 1,dashed]  (v1)
    (x1) edge[line width= \widthedge, out = 20, in = 160, looseness = 1,dashed]  (v2);

    \begin{scope}[on background layer]
        \fill[gray!10] (u1.center) to [out = -20, in = -160, looseness = 1](v1.center) to (x1.center);
        \fill[gray!30] (x1.center) to [out = 20, in = 160, looseness = 1](v2.center) to (v1.center); 
    \end{scope}
    
    \draw (0,-1.5) node{\Large $\bm{S}$};
    \draw (1,-0.5) node{\Large $\bm{S^c}$};
    
\end{tikzpicture}}

\newcommand{\tikzScSSc}[1]{
\begin{tikzpicture}[x=#1 cm, y=#1 cm]
    
    \vertex (a) at (-2,0) [label=90:$p$]{};
    \vertex (b) at (-1,-2) [label=-90:$a$]{};
    \vertex (c) at (0,0) [label=90:$b$]{};
    \vertex (d) at (1,-2) [label=-90:$c$]{};
    \vertex (u) at (2,0) [label=90:$d$]{};
    
    \path
    (a) edge[line width= 1.5]   node[right]{$e^c_1$} (b)
    (b) edge[line width= 1.5]   node[left]{$f^c_1$} node[right]{$e$} (c)
    (c) edge[line width= 1.5]   node[left]{$f$} node[right]{$e^c_2$} (d)
    (d) edge[line width= 1.5]   node[left]{$f^c_2$} (u)
    
    (a) edge[line width= \widthedge, out = 20, in = 160, looseness = 1,dashed]  (c)
    (b) edge[line width= \widthedge, out = -20, in = -160, looseness = 1,dashed]  (d)
    (c) edge[line width= \widthedge, out = 20, in = 160, looseness = 1,dashed]  (u);

    \begin{scope}[on background layer]
        \fill[gray!30] (a.center) to [out = 20, in = 160, looseness = 1](c.center) to (b.center); 
        \fill[gray!10] (b.center) to [out = -20, in = -160, looseness = 1](d.center) to (c.center);
        \fill[gray!30] (c.center) to [out = 20, in = 160, looseness = 1](u.center) to (d.center); 
    \end{scope}
    
    \draw (-1,-0.5) node{\Large $\bm{S^c_1}$};
    \draw (0,-1.5) node{\Large $\bm{S}$};
    \draw (1,-0.5) node{\Large $\bm{S^c_2}$};
    
\end{tikzpicture}}

\begin{document}

\setlength{\mathindent}{2cm} 

\begin{abstract}
A graph $G$ is \emph{$q$-Ramsey} for a $q$-tuple of graphs $(H_1,\ldots,H_q)$ if for every $q$-coloring of the edges of $G$ there exists a monochromatic copy of $H_i$ in color $i$ for some $i\in[q]$. Over the last few decades, researchers have investigated a number of questions related to this notion, aiming to understand the properties of graphs that are $q$-Ramsey for a fixed tuple.
Among the tools developed while studying questions of this type are gadget graphs, called signal senders and determiners, which 
have proven invaluable for building Ramsey graphs with certain properties.
However, until now these gadgets have been shown to exist and used mainly in the two-color setting or in the symmetric multicolor setting, and our knowledge about their existence for multicolor asymmetric tuples is extremely limited. In this paper, we construct such gadgets for any tuple of cliques. We then use these gadgets to generalize three classical theorems in this area to the asymmetric multicolor setting.
\end{abstract}

\maketitle

\section{Introduction}

A graph $G$ is said to be $q$-\emph{Ramsey} for a $q$-tuple of graphs $\cT=(H_1,\ldots,H_q)$, denoted by \mbox{$G\to_q \cT$,} if for every $q$-coloring of the edges of $G$ there exist a color $i\in[q]$ and a monochromatic copy of $H_{i}$ in color $i$. In the special case where $H_i\cong H$ for all $i\in [q]$, we write simply $H$ instead of the tuple $(H,\dots, H)$ and say that  $G$ is $q$-\emph{Ramsey} for $H$. We refer to this setting as the \emph{symmetric case} and simplify  all other notation in a similar way.  The celebrated theorem of Ramsey~\cite{RamseyFormalLogic} establishes the existence of at least one $q$-Ramsey graph for any $q$-tuple of graphs $\cT$. It is then natural to study the collection of $q$-Ramsey graphs for a given $q$-tuple $\cT$. The case where $\cT$ is a tuple of cliques has been of fundamental importance in Ramsey theory and will be the focus of this paper.

A central question in this area asks for the smallest number of vertices in a $q$-Ramsey graph for a given $q$-tuple $\cT$, called the \emph{($q$-color) Ramsey number} of $\cT$ and denoted by $r_q(\cT)$. For tuples of cliques, this question has been studied since the early days of Ramsey theory in the work of \Erdos{} and \Szekeres~\cite{erdos_combinatorial_1935} and \Erdos{}~\cite{erdos1947}. Despite many years of research, even the 2-color Ramsey number $r_2(K_t)$ is far from being understood. The best known bounds for two colors are due to Spencer~\cite{spencer1975ramsey} and Campos, Griffiths, Morris, and Sahasrabudhe~\cite{campos2023exponential} (see also the papers by Sah~\cite{sah2020diagonal} and Conlon~\cite{conlon_upper_bound}). In the multicolor symmetric setting the best known upper bound is obtained using the ideas from~\cite{erdos_combinatorial_1935}, while the lower bound was recently improved in a series of papers by Conlon and Ferber~\cite{conlon_lower_2021}, Wigderson~\cite{wigderson_improved_2021}, and Sawin~\cite{sawin2022improved}. 

This line of research was then generalized to other graph parameters, studying their extremal behavior within the collection of all $q$-Ramsey graphs for a given $\cT$.
One example is the so-called \emph{size-Ramsey number}, denoted by $\hat{r}_q(\cT)$, which represents the minimum number of edges of such a graph. This notion was introduced by \Erdos{}, Faudree, Rousseau, and Schelp~\cite{erdos_size_1978}. Following an argument due to \Chvatal{} (presented in~\cite{erdos_size_1978}), easily generalized to any tuple of cliques $\cT$, we know that $\hat{r}_q(\cT)=\binom{r_q(\cT)}{2}$. Another example appears in the work of Folkman and was
motivated by a question of~\Erdos{} and Hajnal~\cite{research_problems}. \Folkman{}~\cite{Folkman1970Monochromatic} showed that, perhaps surprisingly, for every $t\geq 3$, there exists a graph that is $2$-Ramsey for $K_t$ but contains no copy of $K_{t+1}$, thus establishing the existence of locally sparse Ramsey graphs for a given clique. This result was extended to general graphs and to multiple colors by~\Nesetril{} and \Rodl{}~\cite{Nesetril1976RamseyProperty}.

As any supergraph of a $q$-Ramsey graph for a given $\cT$ is itself $q$-Ramsey for $\cT$, to understand the collection of $q$-Ramsey graphs for $\cT$, it suffices to consider the ones that are minimal with respect to subgraph inclusion. More precisely, we say that a graph $G$ is $q$-\emph{Ramsey-minimal} for $\cT$ if $G$ is $q$-Ramsey for $\cT$ and no  proper subgraph $G'\subsetneq G$ has this property. We denote the set of all such graphs by $\cM_q(\cT)$. In their seminal paper~\cite{burr_graphs_1976}, \BEL{} proved several breakthrough results on the structure and properties of graphs in $\cM_2(K_s,K_t)$. In particular, settling a conjecture of \Nesetril{}, they proved the existence of infinitely many non-isomorphic 2-Ramsey-minimal graphs for $(K_s,K_t)$. They also determined the smallest possible minimum degree, vertex-connectivity, and chromatic number of a graph in $\cM_2(K_s,K_t)$ and showed that there exist $2$-Ramsey-minimal graphs for $(K_s,K_t)$ with arbitrarily large maximum degree, chromatic number, or independence number. 

The proofs of several of these results hinge on the existence of certain gadget graphs, called signal senders, introduced in the same paper. The work of~\BEL{} was extended first by~\Burr{}, Faudree, and Schelp~\cite{burr1977ramseyminimal}, who introduced another closely-related type of gadget, called a determiner. This line of research was continued by \Burr{}, \Nesetril{}, and \Rodl{}~\cite{burr1985useofsenders}, who showed the existence of both types of gadgets for pairs of 3-connected graphs, and finally by~\Rodl{} and~\Siggers{}~\cite{rodl_ramsey_2008}, who proved the existence of signal senders for 3-connected graphs in the symmetric setting for any number of colors. However, the general asymmetric case has been a stumbling point ever since \BEL~\cite{burr_graphs_1976} mentioned that the concepts and questions considered in their paper could be generalized to more than two colors.

In this article, we establish for the first time the existence of signal senders and determiners for arbitrary tuples of cliques. As stating the precise result requires some preparation, we defer the precise statement to the next section (see~\thref{thm:main_existence}). We use these tools to generalize three classical results on the properties of $q$-Ramsey-minimal graphs for tuples of cliques, which we detail now.

\subsection{Ramsey infinite}

In \cite{burr_graphs_1976}, as an immediate consequence of the existence of signal senders, \Burr{}, \Erdos{}, and \Lovasz{} showed that for any integer $t \geq 3$ there are infinitely many non-isomorphic graphs that are $2$-Ramsey-minimal for $K_t$. This result was strengthened first by \Burr{}, \Nesetril{}, and \Rodl{}~\cite{burr1985useofsenders} and more recently by \Rodl{} and \Siggers{}~\cite{rodl_ramsey_2008}, who showed that, for all integers $q\geq2$ and $t \geq 3$ and any sufficiently large $n$, there exist $2^{\Omega(n^2)}$ non-isomorphic graphs on at most $n$ vertices that are $q$-Ramsey-minimal for $K_t$.

In Section~\ref{sec:ramsey_infinite} we prove the following theorem, extending \cite[Theorems 1.1 and 6.1]{rodl_ramsey_2008} to an arbitrary number of colors in the asymmetric case.

\begin{theorem}\thlabel{thm:Ramsey_infinite}
For any $q\geq 2$ and any tuple $\cT = (K_{t_1},\ldots,K_{t_q})$ with $t_1\geq\ldots\geq t_q \geq 3$, there exist constants $c>0$ and $n_0>1$, such that for all $n\geq n_0$, there exist at least $2^{cn^2}$ non-isomorphic graphs on at most $n$ vertices that are $q$-Ramsey-minimal for $\cT$.
\end{theorem}

\subsection{Minimum degree}
\Burr{}, \Erdos{}, and \Lovasz{}~\cite{burr_graphs_1976} also initiated the study of minimum degrees of Ramsey graphs. The parameter $s_q(H_1,\ldots,H_q)$ is defined as the smallest minimum degree among all $q$-Ramsey-minimal graphs for $(H_1,\ldots,H_q)$, that is,
\[ s_q(H_1,\ldots,H_q) = \min \set*{\delta(G) : G\in\cM_q(H_1,\ldots,H_q)}.\]

\Burr{}, \Erdos{}, and \Lovasz{} considered pairs of complete graphs and established that $s_2(K_s,K_t)=(s-1)(t-1)$ for any $s,t\geq 3$. Dealing with symmetric tuples of cliques, \FGLPS{}~\cite{fox2016minimum} showed that $s_q(K_t)$ is quadratic in $q$, up to a polylogarithmic factor, when the size of the clique is fixed. The polylogarithmic factor was settled to be $\Theta(\log q)$ when $t=3$ by Guo and Warnke~\cite{Guo-Warnke20}, following earlier work in~\cite{fox2016minimum}. In the other regime, when the number of colors is fixed, \Han{}, \Rodl{}, and  \Szabo{}~\cite{Han:2018aa} showed that $s_q(K_t)$ is quadratic in the clique size $t$, again up to a polylogarithmic factor. Bounds that are polynomial in both $q$ and $t$ are also known, see~\cite{bamberg2020minimum,bishnoiNewUpperBound2022,fox2016minimum}. 

The symmetric case $s_2(H)$ has been considered for various other classes of graphs. For example, Fox and Lin~\cite{fox_minimum_2007} studied complete bipartite graphs, while \Szabo{}, \Zumstein{}, and \Zurcher{}~\cite{szabo_minimum_2010} extended their work to several other classes of bipartite graphs, including trees and even cycles.

In~\cite{bishnoi2021minimum}, Bishnoi, Clemens, Gupta, Liebenau, and the authors studied the parameter $s_q$ in the asymmetric setting for the first time since the work of \Burr{}, \Erdos{}, and \Lovasz{}, focusing on tuples containing several cliques $K_t$ and several cycles $C_\ell$.

In this paper, we investigate the parameter $s_q$ for general asymmetric tuples of cliques. In Section~\ref{sec:ramsey_mindeg} we prove the following theorem, an extension of the main result in~\cite{bishnoiNewUpperBound2022} for asymmetric tuples of cliques.

\begin{theorem}\thlabel{thm:min_degree}
For any $q\geq 2$ and any tuple $\cT = (K_{t_1},\ldots,K_{t_q})$ with $t_1\geq\ldots\geq t_q \geq 3$, we have
\[ (t_1-1)(t_2-1) \leq s_q(\cT)\leq (8q(t_1-1)\log(t_1-1))^3.\]
\end{theorem}\smallskip

We believe that a similar approach to the one presented in Section~\ref{sec:ramsey_mindeg} could be used to extend the bounds in~\cite{fox2016minimum,Han:2018aa}. As the details would likely be rather technical, for the sake of simplicity we do not pursue this idea further in this paper. 

\subsection{Ramsey equivalence}

The notion of Ramsey equivalence was introduced by~\Szabo{}, \Zumstein, and \Zurcher~\cite{szabo_minimum_2010} (although a question of this flavor was considered in earlier work by Graham, \Luczak{}, \Rodl{}, and \Rucinski{}~\cite{graham2002ramsey}) and arises from the following natural question. Which pairs of graphs $H$ and $H'$ have the same collection of $q$-Ramsey graphs? \FGLPS{}~\cite{fox_what_2014} explored this question when $H$ is a clique. Their result, combined with the earlier work of \Nesetril{} and \Rodl{}~\cite{Nesetril1976RamseyProperty}, shows that there is no connected graph $G\not\cong K_t$ such that $\cM_2(G)=\cM_2(K_t)$. The notion of Ramsey equivalence was further studied for instance in~\cite{axenovich2017,bloom_ramsey_2018,savery2022chromatic}. An asymmetric variant for pairs $(H_1,H_2)$ and $(H'_1,H'_2)$ was recently considered by the first author together with Clemens, Gupta, and Rollin~\cite{boyadzhiyska2024ramsey}. Here we explore the notion of Ramsey equivalence in the multicolor asymmetric setting in more generality, allowing the two tuples to have different lengths.

\begin{definition}\thlabel{def:Ramsey_equivalence}
    Let $q,\ell\geq 2$ be integers, $\cT$ be a $q$-tuple of graphs, and $\cS$ be an $\ell$-tuple of graphs. We say that $\cT$ and $\cS$ are \emph{Ramsey-equivalent} if, for any graph $G$, we have $G\to_{q}\cT$ if and only if $G\to_{\ell}\cS$. In other words, $\cT$ and $\cS$  are Ramsey-equivalent if $\cM_{q}(\cT) = \cM_{\ell}(\cS)$.
\end{definition}

It is easy to check that adding any number of $K_2$'s to a tuple of graphs $\cT$ does not change the collection of Ramsey graphs. We ask which tuples of cliques are Ramsey-equivalent, apart from these degenerate examples. The first step towards addressing this question follows from the work of \Nesetril{} and \Rodl{}~\cite{Nesetril1976RamseyProperty}, who proved that, for every $q\geq 2$ and every graph $H$, there exists a graph $G$ that is $q$-Ramsey for $H$ and has the same clique number as $H$. As a direct consequence, for any tuples of cliques $\cT$ and $\cS$, if the largest cliques in $\cT$ and $\cS$ have different sizes, then $\cT$ and $\cS$ are not Ramsey-equivalent.

Recent developments on the behavior of random graphs can also be used to improve this result. Indeed, the celebrated \Kohayakawa{}-\Kreuter{} conjecture~\cite{kohayakawaThresholdFunctionsAsymmetric1997} seeks an extension of the work of \Rodl{} and \Rucinski{}~\cite{rodl_threshold_1995} to asymmetric tuples. The $1$-statement of this conjecture was proven true by Mousset, Nenadov, and Samotij~\cite{mousset_nenadov_samotij_2020} in full generality, while the $0$-statement was confirmed in the special case of cliques by Marciniszyn, Skokan, Sp\"{o}hel, and Steger~\cite{Marciniszyn_rd_asym_2009} (a full resolution was announced very recently by Christoph, Martinsson, Steiner, and Wigderson~\cite{christoph2024resolution}, following other recent progress by Kuperwasser, Samotij, and Wigderson~\cite{kuperwasser2023kohayakawa} and Bowtell, Hancock, and Hyde~\cite{bowtell2023proof}). Given $\cT=(K_{t_1},\ldots,K_{t_q})$,  these results imply that $n^{-1/m_2(t_1,t_2)}$ is a probability threshold for the binomial random graph\footnote{The binomial random graph is a graph on vertex set $[n]$ in which every possible edge is inserted independently with probability $p$.} $G_{n,p}$ to be $q$-Ramsey for $\cT$, where $m_2(t_1,t_2)$ is a constant depending only on $t_1$ and $t_2$. Using simple properties of $m_2(t_1,t_2)$ (see \cite{Marciniszyn_rd_asym_2009,mousset_nenadov_samotij_2020}), it is easy to deduce that, if $\cS = (K_{s_1},\ldots,K_{s_\ell})$, then $(t_1,t_2)=(s_1,s_2)$ is necessary for $\cT$ and $\cS$ to be Ramsey-equivalent. 

It is then natural to ask whether there exist \emph{any} tuples of cliques that are Ramsey-equivalent. This question was answered in the negative by Graham, \Luczak{}, \Rodl{}, and \Rucinski{}~\cite[Corollary 1.1]{graham2002ramsey}. 
\begin{theorem}[Graham, \Luczak{}, \Rodl{}, and \Rucinski{}~\cite{graham2002ramsey}]\thlabel{thm:Ramsey_equiv}
Let $q,\ell\geq 2$ be integers, and $\cT,\cS$ be two tuples of cliques such that $\cT = (K_{t_1},\ldots,K_{t_q})$ and $\cS = (K_{s_1},\ldots,K_{s_\ell})$ with $t_1\geq\ldots\geq t_q \geq 3$ and $s_1\geq\ldots\geq s_{\ell} \geq 3$. Then $\cT$ and $\cS$ are Ramsey-equivalent if and only if $\cT=\cS$.
\end{theorem}

The original proof relies on a stronger result. Graham, \Luczak{}, \Rodl{}, and \Rucinski{} showed that, given any two tuples of cliques $\cT,\cS$ as in~\thref{thm:Ramsey_equiv}, we have $\cM_q(\cT)\subseteq\cM_{\ell}(\cS)$ if and only if there exists a partition of $\cS$ into $q$ parts $\cS_1,\ldots,\cS_q$, such that for all $i\in[q]$, either $\cS_i$ is empty or $K_{t_i}$ is Ramsey for $\cS_i$.
In Section~\ref{sec:ramsey_equiv}, we use the existence of signal senders and determiners for asymmetric tuples to obtain a new proof of~\thref{thm:Ramsey_equiv}.
\bigskip

\paragraph{\textbf{Organization of the paper.}}
In Section~\ref{sec:prelim}, we introduce key definitions and preliminaries and state our main technical result, \thref{thm:main_existence}. In Section~\ref{sec:using_tools} we prove \thref{thm:Ramsey_infinite,thm:min_degree,thm:Ramsey_equiv}, assuming \thref{thm:main_existence}, which we prove in Section~\ref{sec:existence}.

\section{Preliminaries}\label{sec:prelim}

Much of our notation is standard and similar to~\cite{bishnoi2021minimum}. The \emph{length} of a path in a graph is the number of edges it contains; the \emph{distance} between two edges $uv$ and $xy$ in a graph is the length of a shortest path with one endpoint in $\set*{u,v}$ and one endpoint in $\set*{x,y}$. 

An \emph{$\ell$-graph} is an $\ell$-uniform hypergraph, i.e., a hypergraph in which each hyperedge contains exactly $\ell$ vertices. A \emph{circuit} of length $h$ in a hypergraph is a set of hyperedges $\set*{e_1,\ldots,e_h}$ such that $\size*{\bigcup_{i=1}^he_i}\leq (\ell-1)h$. The \emph{girth} of a hypergraph is the length of its shortest circuit; if no circuit exists, the hypergraph is said to have infinite girth.

A \emph{digraph} is a directed graph $D=(V,E)$, where $V$ is the set of vertices of $D$ and $E$ is a set of ordered pairs of vertices, called \emph{arcs}. For clarity, we often denote  the arc $(i,j)$ by $\vv{ij}$ and write $\vv{ijk}$ for the directed path with arcs $\vv{ij}$ and $\vv{jk}$. A \emph{$2$-cycle} in a digraph $D=(V,E)$ is a pair of arcs $\vv{ij},\vv{ji}\in E$. Given a hypergraph $G$, we let $V(G)$ be its vertex set, and $E(G)$ be its edge set. Given an edge $e\in E(G)$, we write $G-e$ for the spanning subgraph of $G$ obtained by removing the edge $e$. Given hypergraphs $G$ and $H$, we say that $G$ is \emph{$H$-free} if it does not contain $H$ as a not necessarily induced subgraph. 

Given a digraph $D=(V,E)$ and a vertex $u\in V$, we write $N^+_D(u)$ for the out-neighborhood of $u$ in $D$, i.e., $N^+_D(u)=\set*{v\in V : \vv{uv}\in E}$. Similarly, we denote  the in-neighborhood of $u$ in $D$ by $N^-_D(u)$. If $G$ is an undirected graph, we simply write $N_G(u)$ for the neighborhood of $u$ in $G$. 

Unless otherwise specified, $\cT$ will always denote a $q$-tuple of cliques of the form $\cT = (K_{t_1},\dots, K_{t_q})$, where $t_1\geq \dots\geq t_q\geq 3$. Given such a tuple $\cT$ and a subset $X\subseteq [q]$, we write $\cT_{X}$ for the tuple  $(K_{t_i})_{i\in X}$, that is, the $|X|$-tuple consisting of all cliques $K_{t_i}$ such that $i\in X$.  For a subset $X\subseteq[q]$, we write $\overline{X}$ for the set $[q]\setminus X$. 

Unless otherwise specified, the term coloring refers to an edge-coloring. If a coloring of a graph uses at most $q$ colors, then we say that it is a \emph{$q$-coloring}. Typically, the color palette in a $q$-coloring is taken to be the set $[q]=\set*{1,\ldots,q}$. Given a $q$-coloring $\phi$ of a graph $G$ and a subgraph $F\subseteq G$, we will write $\phi_{|F}$ for the $q$-coloring induced by $\phi$ on the edges of $F$. Given a $q$-tuple of graphs $(H_1,\ldots, H_q)$, we say that a $q$-coloring $\phi$ of a graph $G$ is \emph{$(H_1,\ldots, H_q)$-free} if, for all $i\in [q]$, the graph consisting of the edges in $\phi^{-1}(\set{i})$ is $H_i$-free.

Recall that we write $G\to_q\cT$ to denote the fact that $G$ is $q$-Ramsey for a $q$-tuple of graphs $\cT$. Similarly, we write $G\nto_q\cT$ when $G$ is not $q$-Ramsey for $\cT$. When $q$ is clear from context, we will sometimes suppress it from the notation. 

The following result is a simple consequence of the work of \Nesetril{} and \Rodl{}~\cite{Nesetril1976RamseyProperty} discussed in the introduction. It will be handy both in the proof of our main result and in some of our applications. 
\begin{theorem}\thlabel{thm:NR}
Let $q\geq 2$ and $(K_{t_1}, \dots, K_{t_q})$ be a $q$-tuple of cliques such that $t_1\geq \dots\geq t_q\geq 3$. Then there exists a graph $G$ containing no copy of $K_{t_1+1}$ such that $G\to_q (K_{t_1}, \dots, K_{t_q})$.
\end{theorem}
The theorem follows immediately from~\cite[Theorem 1]{Nesetril1976RamseyProperty} by taking a graph $G$ that is $q$-Ramsey for $K_{t_1}$ and contains no copy of $K_{t_1+1}$.

\medskip
We now introduce the two gadget graphs that are the focus of this paper. As discussed above, Burr, Faudree, and Schelp~\cite{burr1977ramseyminimal} introduced so-called determiners for pairs of cliques. The following definition, first used in~\cite{bishnoi2021minimum}, is a suitable generalization of this concept to the multicolor setting.

\begin{definition}[Set-determiner]\thlabel{def:set_determiner}
Let $q\geq 2$ be an integer, $(H_1,\dots, H_q)$ be a $q$-tuple of graphs, and $X\subseteq [q]$ be a non-empty subset of colors. An \emph{$X$-determiner} for $(H_1,\dots, H_q)$  is a graph $R$ with a distinguished edge $e$ satisfying the following properties:
\begin{enumerate}[label=(R\arabic*)]
    \item $R\nrightarrow_q (H_1,\dots, H_q)$.\label{axiom:set_determiner_not_Ramsey}
    \item For any $(H_1,\dots, H_q)$-free coloring $\phi$ of $R$, we have $\phi(e)\in X$.\label{axiom:set_determiner_colored_edge}
    \item For any color $c\in X$, there exists an $(H_1,\dots, H_q)$-free coloring $\phi$ of $R$ such that $\phi(e)=c$.\label{axiom:set_determiner_any_color}
\end{enumerate}
The edge $e$ is referred to as the \emph{signal edge} of $R$. 
\end{definition}

In the special case where $X=\set{c}$ for some color $c$, these gadgets correspond precisely to the \emph{determiners} introduced by \Burr{}, Faudree, and Schelp in~\cite{burr1977ramseyminimal} and shown to exist for pairs of different complete graphs in the same paper. Subsequently, \Burr{}, \Nesetril{}, and \Rodl{}~\cite{burr1985useofsenders} showed that determiners exist for all pairs of non-isomorphic $3$-connected graphs $(G,H)$. More recently, Siggers~\cite{siggers_non-bipartite_2014} showed that determiners exist for some pairs of the form $(C_\ell, H)$, and the authors together with Bishnoi, Clemens, Gupta, and Liebenau~\cite{bishnoi2021minimum} showed that set-determiners exist for some specific tuples of cliques and cycles. Note that colors corresponding to isomorphic graphs in a given tuple $(H_1,\dots, H_q)$ can be permuted, thus implying that an $X$-determiner can only exist for $(H_1,\dots, H_q)$ if $H_i\not\cong H_j$ for any $i\in X, j\not\in X$.
In other words, if $i\in X$, then for any color $j$ such that $H_j\cong H_i$, it is necessary to have $j\in X$ for an $X$-determiner to exist.\smallskip

Set-determiners allow us to pick which set the color of a certain edge should belong to. In order to have control over the specific color pattern we see on a set of edges, we also define the following more sophisticated gadget.

\begin{definition}[Set-sender]\thlabel{def:set_sender}
Let $q\geq 2$ be an integer, $(H_1,\dots, H_q)$ be a $q$-tuple of graphs, and $X\subseteq [q]$ be any subset of colors. A \emph{negative} (respectively \emph{positive}) \emph{$X$-sender} for $(H_1,\dots, H_q)$ is a graph $S$ with distinguished edges $e$ and $f$, satisfying the following properties:
\begin{enumerate}[label=(S\arabic*)]
    \item $S\nrightarrow_q (H_1,\dots, H_q)$.\label{axiom:set_sender_not_Ramsey}
    \item For any $(H_1,\dots, H_q)$-free coloring $\phi$ of $S$, there exist colors $c_1,c_2\in X$, with $c_1\neq c_2$ (respectively $c_1= c_2$), such that $\phi(e)=c_1$ and $\phi(f)=c_2$. \label{axiom:set_sender_colored_edges} 
    \item For any colors $c_1,c_2\in X$, with $c_1\neq c_2$ (respectively $c_1= c_2$), there exists an $(H_1,\dots, H_q)$-free coloring $\phi$ of $S$ with $\phi(e)=c_1$ and $\phi(f)=c_2$.\label{axiom:set_sender_any_color}
\end{enumerate}
The edges $e$ and $f$ are referred to as the \emph{signal edges} of $S$.
\end{definition}

Note that an $X$-sender with $X\subsetneq[q]$ is also an $X$-determiner for the same $q$-tuple with signal edge $e$ (or $f$). In the special case where $q=2$ and $X=[2]$, $X$-senders correspond to the gadgets introduced by \Burr{}, \Erdos{}, and \Lovasz{}~\cite{burr_graphs_1976}, which they called \emph{signal senders}. The existence of signal senders for pairs of complete graphs was established in~\cite{burr_graphs_1976,burr1977ramseyminimal}. As in the case of determiners, this result was generalized to all pairs of 3-connected graphs in~\cite{burr1985useofsenders}. \Rodl{} and Siggers~\cite{rodl_ramsey_2008} and Siggers~\cite{siggers_highly_2008} extended this result to the multicolor setting when $H_i\cong H$ for all $i\in[q]$  and $H$ is either $3$-connected or a cycle. In a later paper, Siggers~\cite{siggers_non-bipartite_2014} constructed signal senders for some pairs of the form $(C_\ell, H)$.\smallskip

These gadgets are used to force specific color patterns on well-chosen sets of edges. In applications, we usually start with some graph $G$ and add set-senders and set-determiners in such a way that, in any $(H_1,\dots,H_q)$-free coloring of the resulting graph, we obtain a particular color pattern on the edges of $G$. More precisely, we will say that we \emph{attach} a set-determiner $R$ with signal edge $e$ to an edge $f$ of $G$ to mean that we create a new graph as the disjoint union of $R$ and $G$, then merge the edges $e$ and $f$ (see Figure~\ref{fig:attachingDet}). In that case, for clarity, we sometimes say that we identify $e$ with $f$. Similarly, we will say that we \emph{join} two edges $a$ and $b$ of $G$ by a set-sender $S$ with signal edges $e$ and $f$ to mean that we create a new graph as the disjoint union of $S$ and $G$, then merge the edges $a$ and $e$ together, and the edges $f$ and $b$ together (see Figure~\ref{fig:joiningSender}). Again, for clarity we usually say that we \emph{identify} $e$ with $a$ and $f$ with $b$.\bigskip

\begin{minipage}{0.49\textwidth}
    \centering
    \tikzAttachDet{2}
    
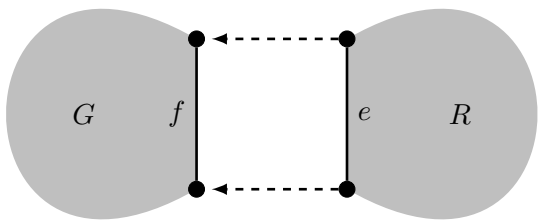
\captionof{figure}{Attaching a set-determiner.}
    \label{fig:attachingDet}
\end{minipage}%
\begin{minipage}{0.49\textwidth}
    \centering
    \tikzJoinSignal{1}
    
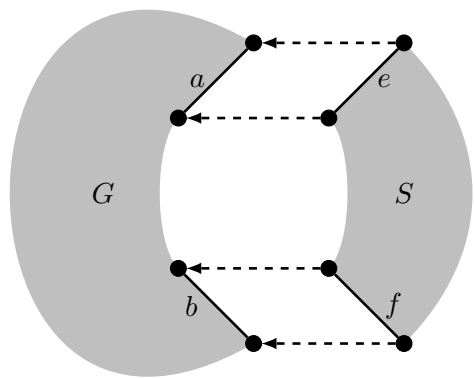
\captionof{figure}{Joining the edges $a$ and $b$ by a set-sender.}
    \label{fig:joiningSender}
\end{minipage}\bigskip

In order for these constructions to be useful, we need to be able to control the new copies of $H_1,\ldots, H_q$ that might be created while adding a set-sender or determiner to some graph. In particular, since we usually use these gadgets as black boxes, we would like to be able to obtain an $(H_1,\dots, H_q)$-free coloring of the entire graph by simply giving each of the building blocks an $(H_1,\dots, H_q)$-free coloring. This motivates the definitions of safe set-senders and safe set-determiners, in line with the concept developed by Siggers in~\cite{siggers_non-bipartite_2014}. 

\begin{definition}[Safe set-determiner]\thlabel{def:safe_set_determiner}
Let $X\subseteq [q]$ be any subset of colors and $R$ be an $X$-determiner for $(H_1,\dots, H_q)$ with signal edge $e$. We say that $R$ is \emph{safe} if Property~\ref{axiom:set_determiner_any_color} can be replaced by the following:
\begin{enumerate}[label=(R\arabic*')]
    \setcounter{enumi}{2}
    \item For any color $c\in X$, there exists an $(H_1,\dots, H_q)$-free coloring $\phi$ of $R$ such that: \label{axiom:safe_set_determiner_any_color}
    \begin{itemize}
        \item $\phi(e)=c$.
        \item If $R$ is attached to any edge of another graph $G$, then a $q$-coloring $\psi$ of $G\cup R$ that extends $\phi$ is $(H_1,\dots, H_q)$-free if and only if $\psi_{\mid G}$ is $(H_1,\dots, H_q)$-free.
    \end{itemize}
\end{enumerate}
\end{definition}

In other words, if a safe $X$-determiner $R$ is attached to some base graph $G$, we can always extend an $(H_1,\dots, H_q)$-free coloring of $G$ to an $(H_1,\dots, H_q)$-free coloring of the entire graph $G\cup R$, provided that the signal edge of $R$ receives one of the colors from the set $X$. A safe set-sender is defined in a similar way.  

\begin{definition}[Safe set-sender]\thlabel{def:safe_set_sender}
Let $X\subseteq [q]$ be any subset of colors and $S$ be a negative (respectively positive) $X$-sender for $(H_1,\dots, H_q)$ with signal edges $e,f$. Then $S$ is said to be \emph{safe} if Property~\ref{axiom:set_sender_any_color} can be replaced by the following:
\begin{enumerate}[label=(S\arabic*')]
    \setcounter{enumi}{2}
    \item For any colors $c_1,c_2\in X$ with $c_1\neq c_2$ (respectively $c_1= c_2$), there exists an $(H_1,\dots, H_q)$-free coloring $\phi$ of $S$ such that:\label{axiom:safe_set_sender_any_color}
    \begin{itemize}
        \item $\phi(e)=c_1$ and $\phi(f)=c_2$.
        \item If $S$ joins any two edges of another graph $G$, then a $q$-coloring $\psi$ of $G\cup S$ that extends $\phi$ is $(H_1,\dots, H_q)$-free if and only if $\psi_{\mid G}$ is $(H_1,\dots, H_q)$-free.
    \end{itemize}
\end{enumerate}
\end{definition}

We now prove that set-determiners for tuples of cliques are always safe, while set-senders for such tuples are safe provided that the two signal edges are sufficiently far apart. The proofs are similar to those provided for example in~\cite{bishnoi2021minimum}, so we only sketch them.

\begin{lemma}\thlabel{thm:all_det_are_safe}
For any $q\geq 2$ and any tuple $\cT = (K_{t_1},\ldots,K_{t_q})$ with $t_1\geq\ldots\geq t_q \geq 3$, any set-determiner for $\cT$ is safe.
\end{lemma}

\begin{proof}
Let $R$ be an $X$-determiner for $\cT$ for some subset $\emptyset\subsetneq X\subsetneq [q]$ and $e$ be its signal edge. Suppose $R$ is attached to some edge of a graph $G$. Notice that any clique of size at least three is fully contained either in $R$ or in $G$. This fact together with Property~\ref{axiom:set_determiner_any_color} proves the claim.
\end{proof}

\begin{lemma}\thlabel{thm:all_senders_are_safe}
For any $q\geq 2$ and any tuple $\cT = (K_{t_1},\ldots,K_{t_q})$ with $t_1\geq\ldots\geq t_q \geq 3$, if $S$ is a set-sender for $\cT$ with signal edges $e$ and $f$ such that the distance between $e$ and $f$ is at least three, then $S$ is safe.
\end{lemma}

\begin{proof}
The proof is very similar to that of the previous lemma. Let $S$ be an $X$-sender for $\cT$ for some subset $\emptyset\subsetneq X$ and $e$ and $f$ be its signal edges. Suppose $S$ joins two edges of a graph $G$. Again, since the distance between $e$ and $f$ is at least three, any clique of size at least three is fully contained either in $S$ or in $G$. This fact together with Property~\ref{axiom:set_sender_any_color} proves the claim.
\end{proof}

We finally introduce the concept of a \emph{distinguishable} tuple.
Informally, we call a tuple $\cT$ distinguishable if all ``useful'' set-determiners and set-senders for $\cT$ exist and can be chosen to be safe.

\begin{definition}[Distinguishable]\thlabel{def:distinguishable}
Let $q\geq 2$ and $\cT=(H_1,\ldots,H_q)$ be a $q$-tuple of graphs. For any graph $F$, define $X_F\subseteq[q]$ to be the set $X_F=\set*{i\in[q] : H_i\cong F}$. We say that $(H_1,\ldots,H_q)$ is \emph{distinguishable} if, for any graph $F$ such that $X_F\neq\emptyset$, the following hold: 
\begin{enumerate}
    \item There exists a safe $X_F$-determiner for $\cT$.
    \item If $\size{X_F}>1$, then there exist safe positive and negative $X_F$-senders for $\cT$.
\end{enumerate}
\end{definition}
When working with tuples of cliques, we will often simplify notation and write $X_i$ instead of $X_{K_{t_i}}$. For example, by the above discussion we know that symmetric tuples $(H,\dots, H)$, where $H$ is 3-connected or isomorphic to a cycle, are distinguishable, and that pairs of non-isomorphic 3-connected graphs are distinguishable. 
As explained above, if $X\cap X_F$ is a nontrivial subset of $X_F$ for some graph $F$, then no $X$-sender or $X$-determiner can exist (as we can simply permute the colors in $X_F$). Therefore, a tuple of graphs is distinguishable if the ``most restrictive'' possible gadgets exist.\smallskip

We are ready to state our main result, establishing that any tuple consisting of nontrivial cliques is distinguishable. 

\begin{theorem}\thlabel{thm:main_existence}
Let $q\geq2$ and $t_1\geq \ldots \geq t_q\geq 3$. The tuple $(K_{t_1}, \dots, K_{t_q})$ is distinguishable.
\end{theorem}\medskip

\section{Applications of set-senders and set-determiners}\label{sec:using_tools}

Before proving our main result (\thref{thm:main_existence}), we demonstrate its utility by proving~\thref{thm:Ramsey_infinite,thm:min_degree,thm:Ramsey_equiv} as applications. We assume that~\thref{thm:main_existence} holds throughout the entire section. The results in this section mostly generalize known theorems from the symmetric setting, using set-senders and set-determiners instead of classical signal senders and determiners. For simplicity, when this is the case, we present only sketches of the proofs here. Readers familiar with previous work on this topic are free to skip this section.

\subsection{Ramsey infinite}\label{sec:ramsey_infinite}

In this section we prove \thref{thm:Ramsey_infinite} assuming~\thref{thm:main_existence}. This strengthens \cite[Theorems 1.1 and 6.1]{rodl_ramsey_2008} and extends these results to an arbitrary number of colors in the asymmetric case.

\begin{prop}\thlabel{thm:Ramsey_infinite_q}
Let $q\geq 2$ and $\cT=(K_{t_1},\ldots,K_{t_q})$ be a $q$-tuple of cliques such that ${t_1\geq\ldots\geq t_q\geq 3}$. Assume that any $q$-tuple of cliques is distinguishable. Then there exist constants $c=c(\cT)>0$ and $n_0=n_0(\cT)>1$ such that, for all $n\geq n_0$, there exist at least $2^{cn^2}$ non-isomorphic graphs on at most $n$ vertices that are $q$-Ramsey-minimal for $\cT$.
\end{prop}

The proof of this proposition relies on the following lemma.

\begin{lemma}\thlabel{lemma:Ramsey_infinite_induction}
Let $1\leq i<q$ and $t_1,\ldots,t_q$ be integers such that $t_1 = \ldots = t_i\geq t_{i+1}\geq\ldots\geq t_q\geq 3$. Assume that any $q$-tuple of cliques is distinguishable. Then there exists a constant $s=s(q,t_1,\ldots,t_q)$ such that for any graph $F\in\cM_q(K_{t_1},\ldots,K_{t_q})$, there exists a graph $G$ with the following properties:
\begin{enumerate}
    \item $G$ has at most $s\size{V(F)}$ vertices.
    \item $G$ is $q$-Ramsey for $(K_{t_1+1},\ldots,K_{t_i+1},K_{t_{i+1}},\ldots,K_{t_q})$.
    \item Any subgraph $G'\subseteq G$ such that $G'\in\cM_q(K_{t_1+1},\ldots,K_{t_i+1},K_{t_{i+1}},\ldots,K_{t_q})$ contains $F$ as a subgraph.
\end{enumerate}
\end{lemma}

Given $t_1=\ldots =t_i\geq t_{i+1}\geq\ldots\geq t_q$ and $\cT^+=(K_{t_1+1},\ldots,K_{t_i+1},K_{t_{i+1}},\ldots,K_{t_q})$, the proof of this lemma is a direct generalization of the work of \Rodl{} and \Siggers{}, \cite[Lemma 6.2]{rodl_ramsey_2008}, using $[i]$-senders for $\cT^+$ when $i\in\{2,\ldots,q-1\}$ and $\{1\}$-determiners for $\cT^+$ otherwise.

The proof of \thref{thm:Ramsey_infinite_q} then follows by repeated induction using \thref{lemma:Ramsey_infinite_induction}, and is again a generalization of \cite[Theorem 6.1]{rodl_ramsey_2008}. The process can be summarized as follows. The base case of the induction deals with the symmetric tuple $(K_{t_q},\ldots,K_{t_q})$, for which the theorem is true by \cite[Theorem  1.1]{rodl_ramsey_2008}. We then let $i_1$ be the largest index such that $t_{i_1}>t_q$ in $\cT$. By induction, repeated $(t_{i_1+1}-t_{i_1})$ times, applying~\thref{lemma:Ramsey_infinite_induction} we obtain \thref{thm:Ramsey_infinite_q} for the tuple $(K_{t_{i_1}},\ldots,K_{t_{i_1}},K_{t_{{i_1}+1}},\ldots,K_{t_q})$, with ${i_1}$ cliques of size $t_{i_1}$ and $t_{{i_1}+1}=\ldots=t_q$. We then repeat the induction, choosing ${i_2}$ to be the largest index such that $t_{i_2}>t_{i_1}$, yielding \thref{thm:Ramsey_infinite_q} for the tuple $(K_{t_{i_2}},\ldots,K_{t_{i_2}},K_{t_{{i_2}+1}},\ldots,K_{t_q})$. We repeat the process for each new clique size, until we obtain \thref{thm:Ramsey_infinite_q} for the tuple $(K_{t_1},\ldots,K_{t_q})$. For example, to prove \thref{thm:Ramsey_infinite_q} for the tuple $(K_5,K_5,K_4,K_3,K_3)$, the repeated induction works with the following tuples:

\[ \underbrace{(K_3,K_3,K_3,K_3,K_3)}_{\text{Base case}}
\xRightarrow[\text{Lemma }\ref{lemma:Ramsey_infinite_induction}]{}
\underbrace{(K_4,K_4,K_4,K_3,K_3)}_{i_1=3,\ t_3=4}
\xRightarrow[\text{Lemma }\ref{lemma:Ramsey_infinite_induction}]{}
\underbrace{(K_5,K_5,K_4,K_3,K_3)}_{i_2=2,\ t_2=5}.\]

The final counting argument is identical to the one used by \Rodl{} and \Siggers{} in \cite[Theorem 6.1]{rodl_ramsey_2008}.

\subsection{Minimum degree}\label{sec:ramsey_mindeg}

In this section we prove \thref{thm:min_degree} assuming \thref{thm:main_existence}. We start by generalizing the packing parameter defined in~\cite{fox2016minimum}. A {\em color pattern} on vertex set $V$ is a collection of edge-disjoint graphs $G_1,\ldots,G_m$ on the same vertex set $V$.

\begin{definition}[Packing parameter]\thlabel{def:packing}
For given integers $q\geq 2$ and $t_1\geq\ldots\geq t_q\geq 2$, we define the packing parameter $P_q(t_1,\ldots,t_q)$ to be the smallest integer $n$ such that there exists a color pattern $G_1,\ldots,G_q$ on vertex set $[n]$ satisfying the following properties:
\begin{enumerate}
    \item $G_i$ is $K_{t_i+1}$-free for every $i\in[q]$. \label{item:Packing_P1}
    \item For every vertex-coloring $\lambda:[n]\to[q]$, there exists a color $i\in[q]$ such that $G_i$ contains a copy of $K_{t_i}$ on the vertices of color $i$.
    \label{item:Packing_P2}
\end{enumerate}
\end{definition}

For $t_1=\ldots=t_q$, this parameter was introduced in~\cite{fox2016minimum}, where it was shown that, for all $q\geq 2$ and $t\geq 3$, we have $s_{q}(K_t) = P_{q}(t-1,\ldots,t-1)$. We generalize this result to arbitrary tuples of cliques. 

\begin{prop}\thlabel{lem:packing_equivalence}
Let $q\geq 2$ and $\cT=(K_{t_1},\ldots,K_{t_q})$ be a $q$-tuple of cliques such that $t_1\geq\ldots\geq t_q\geq 3$. If $\cT$ is distinguishable, then
\[ s_{q}(\cT) = P_{q}(t_1-1,\ldots,t_q-1).\]
\end{prop}

The proof of the lower bound $s_q(\cT) \geq P_{q}(t_1-1,\ldots,t_q-1)$ is identical to the symmetric case (see~\cite[Lemma 2.1]{fox2016minimum}). The proof of the upper bound $s_q(\cT) \leq P_{q}(t_1-1,\ldots,t_q-1)$ is a direct generalization of \cite[Theorem 2.3]{fox2016minimum}, using the set-senders and set-determiners guaranteed by~\thref{thm:main_existence} instead of signal senders.

We are now ready to prove \thref{thm:min_degree}. The lower bound is a direct argument related to the monotonicity of the parameter $s_q$ in $q$. 

\begin{lemma}\thlabel{lem:packing_lower_bound}
Let $q\geq 2$ and $t_1\geq\ldots\geq t_q\geq 3$. Then
\[ t_1t_2\leq P_{q}(t_1,\ldots,t_q).\]
\end{lemma}

\begin{proof}
Let $n=P_q(t_1,\ldots,t_q)$, and $G_1,\ldots,G_q$ be a color pattern as guaranteed by \thref{def:packing}. Then for every vertex-coloring $\lambda:[n]\to[q]$, there exists a color $i\in[q]$ such that $G_i$ contains a copy of $K_{t_i}$ on the vertices of color $i$. This is true in particular for any such coloring using only the two colors $\set{1,2}$. Therefore $G_1,G_2$ is a color pattern satisfying:
\begin{enumerate}
    \item for every $i\in\{1,2\}$, the graph $G_i$ is $K_{t_i+1}$-free, and
    \item for every vertex-coloring $\lambda:[n]\to\{1,2\}$, there exists a color $i\in\{1,2\}$ such that $G_i$ contains a copy of $K_{t_i}$ on the vertices of color $i$.
\end{enumerate}
Hence $n\geq P_2(t_1,t_2)$. The work of~\BEL{}~\cite{burr_graphs_1976} establishes that $P_2(t_1,t_2)=t_1t_2$ and we conclude that $n\geq t_1t_2$.
\end{proof}

The upper bound of \thref{thm:min_degree} requires asymmetric generalizations of some known lemmas. Following a methodology initially developed by \Dudek{} and \Rodl{}~\cite{Dudek:2011aa} and then by ~\FGLPS{}~\cite{fox2016minimum}, the following lemma was proven by Bamberg, Bishnoi, and the second author~\cite{bamberg2020minimum} in the special case $t_1=\dots=t_q$. While the original result is stated in the language of finite geometry, we translate it here to the language of hypergraphs to avoid introducing a whole new set of definitions.

\begin{lemma}\thlabel{lemma:asymmetric_BBL}
Let $q\geq 2$ and $t_1\geq\ldots\geq t_q\geq 2$. Let $s, k$ be positive integers. Suppose there exists a family $\cH_1,\dots,\cH_q$ of edge-disjoint $s$-uniform hypergraphs on the same vertex set $\cV$, each of which is $k$-regular and has girth at least four. Suppose further that the hypergraph on $\cV$ with edge set $\bigcup_{i=1}^qE(\cH_i)$ has girth at least three. If $s > 3qt_1\ln t_1$ and $k  > 3t_1(1+\ln q)$, then $P_q(t_1,\ldots,t_q) \leq \size{\cV}$. 
\end{lemma}

For simplicity, we present only a sketch of the proof here, explaining why the original lemma can be extended to the asymmetric setting.

\begin{proof}
The main ideas behind the construction when $t_1=\dots=t_q = t$ are as follows. Let $n=\size{\cV}$ and without loss of generality assume that $\cV=[n]$. For every hyperedge $\mathcal{E}$ in every hypergraph $\cH_i$, take a uniformly random equipartition\footnote{A partition where the sizes of any two parts differ by at most one.} of $\mathcal{E}$ into $t$ parts, each choice of partition being made independently. For every $i\in[q]$, let $G_i$ be the graph on vertex set $[n]$ obtained by adding the following edges: for each hyperedge $\mathcal{E}\in\mathcal{H}_i$, add a complete $t$-partite graph on the vertices of $\mathcal{E}$, respecting the random partition of $\mathcal{E}$ chosen earlier. Since the girth condition on $\cH_i$ implies that any pair of hyperedges share at most one vertex, $G_i$ is then the union of disjoint \Turan{} graphs. Each $G_i$ is $K_{t+1}$-free, and one can show that with the appropriate condition on $s$ and $k$, for every $i\in [q]$, with positive probability every set of at least $n/q$ vertices contains a copy of $K_t$ in $G_i$.
This implies (by the pigeonhole principle) that, for every $q$-vertex coloring of $[n]$, there exists a color $i\in [q]$ such that $G_i$ contains a copy of $K_{t}$ all of whose vertices have color $i$. Hence $n\geq P_q(t,\ldots,t)$. We refer to~\cite[Section 3]{bamberg2020minimum} for the details.

In the asymmetric setting, the conditions $s > 3qt_1\ln t_1$ and $k  > 3t_1(1+\ln q)$ ensure that for all $i\in[q]$, we have $s > 3qt_i\ln t_i$ and $k > 3t_i(1+\ln q)$. As all choices of random partitions are made independently, we apply the same strategy as in the symmetric case, where the partitions of the hyperedges in $\cH_i$ now contain $t_i$ parts (instead of the constant $t$ parts in the symmetric setting). This yields a color pattern $G_1,\ldots,G_q$ such that, for all $i\in[q]$, the graph $G_i$ is $K_{t_i+1}$-free and every set of at least $n/q$ vertices contains a copy of $K_{t_i}$ in $G_i$. By the pigeonhole principle, for any $q$-vertex coloring of $[n]$, there exists a color $i\in[q]$ used on at least $n/q$ vertices. Thus, $G_i$ contains a copy of $K_{t_i}$ on the vertices of color $i$, and $n\geq P_q(t_1,\ldots,t_q)$.
\end{proof}

Using a construction from \cite{fox2016minimum}, Bishnoi and the second author proved the following lemma in~\cite{bishnoiNewUpperBound2022}.

\begin{lemma}[Bishnoi and Lesgourgues~\cite{bishnoiNewUpperBound2022}]\thlabel{lemma:exist_PLS}
Let $p$ be any prime power. There exists a family $\cH_1,\ldots,\cH_{p-1}$ of edge-disjoint $p$-uniform hypergraphs on the same vertex set $\cV$ of size $p^3$, each of which is $(p-1)$-regular and has girth at least four, and such that the hypergraph on $\cV$ with edge set $\bigcup_{i=1}^{p-1}E(\cH_i)$ has girth at least three.
\end{lemma}

\thref{lemma:asymmetric_BBL,lemma:exist_PLS} now imply the upper bound of \thref{thm:min_degree}.

\begin{cor}\thlabel{cor:packing_upper_bound}
Let $q\geq 2$ and $t_1\geq\ldots\geq t_q\geq 2$. Then
\[P_{q}(t_1,\ldots,t_q) \leq (8qt_1\log t_1)^3.\]
\end{cor}

\begin{proof}
Let $q\geq 2$ and $t_1\geq\ldots\geq t_q\geq 2$, and let $p$ be the smallest prime such that $p\geq 4t_1q\log t_1$. By Bertrand's postulate, $p\leq 8t_1q\log t_1$. By~\thref{lemma:exist_PLS}, there exists a family of $q<p$ edge-disjoint $p$-uniform hypergraphs on the same vertex set $\cV$, each of which is $(p-1)$-regular and has girth at least four, and such that the hypergraph on $\cV$ with edge set $\bigcup_{i=1}^qE(\cH_i)$ has girth at least three. Note that, for $t_1\geq 2$ and $q \geq 2$, we have $p-1 \geq 3qt_1\ln t_1$ and $p-2 \geq 3t_1(1+\ln q)$. By~\thref{lemma:asymmetric_BBL}, $P_q(t_1,\ldots,t_q)\leq \size{\cV}$, and then $\size{\cV}= p^3$ yields the desired bound.
\end{proof}

It follows from \thref{lem:packing_equivalence}, \thref{lem:packing_lower_bound}, and \thref{cor:packing_upper_bound}  that for any distinguishable tuple of cliques $\cT=(K_{t_1},\ldots,K_{t_q})$ with $t_1\geq\ldots\geq t_q\geq 3$, we have \[(t_1-1)(t_2-1)\leq s_q(\cT)\leq (8q(t_1-1)\log (t_1-1))^3.\] Then \thref{thm:min_degree} follows immediately from \thref{thm:main_existence}.\medskip

\subsection{Ramsey equivalence}\label{sec:ramsey_equiv}
Recall that \thref{thm:Ramsey_equiv} states that any two distinct tuples of cliques are not Ramsey-equivalent. Assuming \thref{thm:main_existence}, this result is a direct consequence of the following lemma. Note that, to avoid convoluted notation, in this section we will assume the cliques in a given tuple to be ordered in a \textbf{nondecreasing} fashion.

\begin{prop}\thlabel{thm:Ramsey_equivalence_q}
For any $q,\ell \geq 2$, let $\cT=(K_{t_1},\ldots,K_{t_{q}})$ and $\cS=(K_{s_1},\ldots,K_{s_{\ell}})$ be two tuples of cliques, where $3\leq t_1\leq\ldots\leq t_{q}$ and $3\leq s_1\leq\ldots\leq s_{\ell}$. Suppose that $\cT$ and $\cS$ are both distinguishable. Then $\cT$ and $\cS$ are Ramsey-equivalent if and only if $\cT=\cS$.
\end{prop}

\begin{proof}

We reason by induction on $q+\ell$. Note that the statement is trivially true for $q+\ell=2$.

Assume now that the statement is true up to $q+\ell-1$. Consider two tuples $\cT=(K_{t_1},\ldots,K_{t_q})$ and $\cS=(K_{s_1},\ldots,K_{s_{\ell}})$ with $3\leq t_1\leq\ldots\leq t_q$ and $3\leq s_1\leq\ldots\leq s_\ell$.

Trivially if $\cT=\cS$ then $\cT$ and $\cS$ are Ramsey-equivalent. Assume now that $\cT\neq\cS$. Our goal is to exhibit a graph $G$ such that $G$ is $q$-Ramsey for $\cT$ but not $\ell$-Ramsey for $\cS$ (or vice versa).

Assume first that $t_i=s_i$ for all $i\in [\min\set{q,\ell}]$. If $q<\ell$, let $G\in\cM_q(\cT)$ and let $e \in E(G)$. By minimality, we can find  a $\cT$-free $q$-coloring $\phi$ of $G-e$, and we extend $\phi$ to $e$ by setting $\phi(e) = q+1\leq \ell$. As $t_i=s_i$ for all $i\in [q]$ and $s_{q+1}\geq 3$, this is an $\cS$-free $\ell$-coloring of $G$, i.e., $G\not\in\cM_\ell(\cS)$, and hence $\cT$ and $\cS$ are not Ramsey-equivalent. The case where $q>\ell$ is handled similarly.

We can therefore assume that there exists an index $k\in [q]$ such that $t_k\neq s_k$. Let $i$ be the smallest such index, and without loss of generality assume that $t_i<s_i$. Note that we therefore have 
\[ (t_1=s_1) \leq (t_2=s_2)\leq \ldots \leq (t_i<s_i).\]

As $\cT$ is distinguishable, let $R$ be a safe $X$-determiner for $\cT$, with signal edge $e$, where $X = \set*{j\in[q] :  t_j=t_i}$, that is, $X$ is the set of all colors associated with cliques of size $t_i$ in $\cT$.
By definition, $R\nto_{q}\cT$ and therefore we can assume that $R\nto_{\ell}\cS$, as otherwise we are done. Recall that $\cT_{X}$ denotes the tuple $(K_{t_i})_{i\in X}$. We consider two cases.\smallskip

\noindent\textbf{Case 1}: Assume that there exists an $\cS$-free $\ell$-coloring $\phi$ of $R$ such that $\phi(e)=j\geq i$. Let $F$ be a $K_{t_i+1}$-free graph such that $F\in\cM_{\size{X}}(\cT_{X})$ (as guaranteed by \thref{thm:NR}). Note that $F$ is also $K_{s_j}$-free, as $t_i<s_i\leq s_j$. Attach a copy of $R$ to each edge of $F$. Call the resulting graph $G$. As $R$ is an $X$-determiner for $\cT$, we have $G\to_{q}\cT$. On the other hand, we claim that $G\nto_{\ell}\cS$. Indeed, apply $\phi$ to each copy of $R$ so that $F$ is monochromatic in color $j$. As $F$ is $K_{s_j}$-free and every clique on at least three vertices is fully contained either in $F$ or in a copy of $R$, this is an $\cS$-free $\ell$-coloring of $G$. 

\smallskip\label{item:equivalence_higher_color}

\noindent\textbf{Case 2}: We can assume now that for any $\cS$-free $\ell$-coloring $\phi$ of $R$, we have $\phi(e)<i$, that is $R$ is a $Y$-determiner for $\cS$, for some $Y\subseteq\set{1,\ldots,i-1}$. As $\size{X}+\size{Y}<q+l$, by induction we know that $\cT_{X}$ and $\cS_{Y}$ are not Ramsey-equivalent. Again we consider two possibilities.

\hspace{1cm}\textbf{Case 2.1}: Assume that there exists a graph $F$ such that $F\to \cT_{X}$ and $F\nto\cS_{Y}$. Then attach a copy of $R$ to each edge of $F$ and call the resulting graph $G$. Because $R$ is an $X$-determiner for $\cT$ it follows that $G\to_{q}\cT$. However, since $F\nto\cS_{Y}$, there exists  an $\cS$-free $\size{Y}$-coloring $\psi$ of $F$, using colors in $Y$. Because $R$ is also a $Y$-determiner for $\cS$, the coloring $\psi$ can be extended to all copies of $R$ to an $\cS$-free coloring of $G$, implying that $G\nto_{\ell}\cS$.

\hspace{1cm}\textbf{Case 2.2}: Assume that there exists a graph $F$ such that $F\nto\cT_{X}$ and $F\to\cS_{Y}$. Attach a copy of $R$ to each edge of $F$ and call the resulting graph $G$.
The argument is identical to the previous one, using the fact that $R$ is a $Y$-determiner for $\cS$ to deduce that $G\to_{\ell}\cS$, and that $R$ is an $X$-determiner for $\cT$ to conclude that $G\nto_{q}\cT$.
\end{proof}

\section{Proof of~\texorpdfstring{\thref{thm:main_existence}}{Theorem 2.9}}\label{sec:existence}

In this section we prove~\thref{thm:main_existence}, that is, we construct set-determiners and set-senders for any tuple of cliques. The proof will be by strong induction on the number of colors $q$. In Section~\ref{sec:senders_to_determiners}, we show that the existence of any one set-sender or set-determiner, combined with the induction hypothesis, yields the existence of all required gadgets. Section~\ref{sec:basic_construction} introduces a standard construction due to \Rodl{} and \Siggers{}~\cite[Section 7]{rodl_ramsey_2008}. In the symmetric setting, this construction yields a negative $[q]$-sender. However, for asymmetric tuples this construction does not necessarily produce a useful gadget immediately. The issue is that, since colors cannot be permuted as in the symmetric setting, we might not be able to guarantee that any ``permissible'' pair of colors can appear on the signal edges, that is, Property~\ref{axiom:set_sender_any_color} might fail to hold. We overcome this challenge in Section~\ref{sec:constructions_new}, using a series of constructions, starting with the basic one from Section~\ref{sec:basic_construction}, and with each new construction exhibiting ``nicer'' properties and bringing us closer to our goal.

Before we begin, we recall some notation, conventions, and facts that will be useful throughout. Unless otherwise specified, $\cT$ will always denote a $q$-tuple of cliques $(K_{t_1},\dots, K_{t_q})$, where ${t_1\geq \dots \geq t_q\geq 3}$. We will often need to restrict our tuple. For a subset $X\subseteq [q]$, we write $\cT_{X}$ for the $|X|$-tuple containing only the cliques in $\cT$ whose indices are in $X$.
When we refer to a $q$-coloring, we will normally assume that our color palette is the set $[q]$. However, when we work with a larger $q$-tuple $\cT$ and we talk about a $\cT_{X}$-free coloring of some graph, for simplicity our convention will be that we color using the color palette $X$ (as opposed to $[|X|]$). For example, if $\cT=(K_a,K_b,K_c)$ and $X=\set*{1,3}$, we have $\cT_{X} = (K_a,K_c)$ and a $\cT_{X}$-free coloring of some graph is a $2$-coloring with colors $\set*{1,3}$ with no $K_a$ in color $1$ and no $K_c$ in color $3$.

\subsection{From set-senders to set-determiners and back}\label{sec:senders_to_determiners}

The proof of \thref{thm:main_existence} will be by strong induction on the number of colors $q$. We will assume that any $q'$-tuple  of cliques with $q'< q$ is distinguishable (cf.\ \thref{def:distinguishable}). In this subsection, we show that in order to establish that a $q$-tuple is distinguishable, it suffices to be able to construct \emph{any one} set-sender or set-determiner, that is, the existence of any set-sender or set-determiner, combined with the induction hypothesis yields the existence of all required gadgets.

\begin{prop}\thlabel{thm:one_det_is_enough}
Let $q\geq 3$ and $\cT=(K_{t_1},\ldots,K_{t_q})$ with $t_1\geq \ldots \geq t_q\geq 3$. Assume that any $q'$-tuple of cliques with $q'<q$ is distinguishable and that, for some non-empty set of colors $X\subsetneq [q]$, there exists a safe $X$-determiner for $\cT$. Then $\cT$ is distinguishable.
\end{prop}

\begin{proof}
Let $R$ be a safe $X$-determiner for $\cT$ for some $X \subsetneq [q]$. For each $i\in[q]$, define the set $X_i = \set*{j\in[q] :t_i=t_j}$. Notice that, by permuting the colors in the set $X_i$ for each $i\in[q]$, we can conclude that, for each $i\in[q]$, either $X_i\subseteq X$ or $X_i\cap X = \emptyset$.

First we build an $\overline{X}$-determiner $R^c$ for $\cT$, where we recall that $\overline{X}$ denotes the set $[q]\setminus X$. To do so, let $G$ be an $|X|$-Ramsey-minimal graph for $\cT_{X}$ and let $e$ be an arbitrary edge of $G$. Now we obtain $R^c$ from $G$ by attaching a copy of $R$ to every edge of $G$ except for $e$. Then using the properties of $R$ and its safeness, it is not difficult to verify that $R^c$ is indeed an $\overline{X}$-determiner with signal edge $e$. Moreover, by~\thref{thm:all_det_are_safe}, $R^c$  is safe.

Now, by the induction hypothesis, we know that $\cT_{X}$ and $\cT_{\overline{X}}$ are both distinguishable. Let $i\in[q]$ and without loss of generality assume that $X_i\subseteq X$ (the case $X_i\subseteq \overline{X}$ is similar using $R^c$ instead of $R$).  Now, let $R_i$ be a safe $X_i$-determiner and if $|X_i|>1$ let $S_i$ be a safe positive or negative $X_i$-sender for $\cT_{X}$ (with color palette $X$). By attaching a copy of $R$ to every edge of $R_i$ or $S_i$, we then obtain the required $X_i$-determiner and positive/negative $X_i$-sender for $\cT$ with the same signal edges as $R_i$ and $S_i$, respectively. The safeness of these gadgets also follows easily from the safeness of each building block.\end{proof}

\begin{prop}\thlabel{thm:one_sender_is_enough}
Let $q\geq 3$ and $\cT=(K_{t_1},\ldots,K_{t_q})$ with $t_1\geq \ldots \geq t_q\geq 3$. Assume that any $q'$-tuple of cliques with $q'<q$ is distinguishable and that, for some set of colors $\emptyset\subsetneq X\subseteq [q]$, there exists a safe $X$-sender for $\cT$. Then $\cT$ is distinguishable.
\end{prop}

\begin{proof}
Let $S$ be a safe $X$-sender for $\cT$. If $X\subsetneq [q]$, we can treat $S$ as an $X$-determiner whose signal edge is either of the signal edges of $S$ and apply~\thref{thm:one_det_is_enough}. 

Assume then that $X = [q]$. If $t_1 = t_q$, following the work of \Rodl{} and \Siggers{}~\cite{rodl_ramsey_2008}, we know that positive and negative $[q]$-senders for $\cT$ exist with arbitrarily large distance between the signal edges. By~\thref{thm:all_senders_are_safe}, these senders are safe. So we may further assume that $t_1>t_q$. 

First suppose $S$ is a positive $X$-sender for $\cT$. Let $G$ be a copy of $K_{t_q}$ and $e$ be an edge disjoint from $G$. Join $e$ to every edge of $G$ by a copy of $S$ to obtain the graph $\widetilde{G}$. It is then not difficult to check that $\widetilde{G}$ is a $Y$-determiner for $\cT$ with signal edge $e$ for some set $Y$ containing 1 but not containing $q$. Indeed, if we assign color 1 to $e$ and each edge of $G$, then we can extend this to a $\cT$-free coloring of all of $\widetilde{G}$ (by the safeness of each signal sender). In addition, $e$ cannot receive color $q$ in a $\cT$-free coloring of $\widetilde{G}$, as in that case all edges of $G$ would also receive color $q$, a contradiction. Hence, we can apply~\thref{thm:one_det_is_enough} to complete the proof.  

Finally, suppose $S$ is a negative  $X$-sender for  $\cT$. Let $e_1,\dots, e_q, e_{q+1}$ be a matching. For any pair $(i,j)$ with $1\leq i< j\leq q+1$ except the pair $(i,j)=(q,q+1)$, join the edges $e_i$ and $e_j$ by a copy of $S$, calling the newly constructed graph $\widetilde{S}$. It is then not difficult to verify that $\widetilde{S}$ is a safe positive $X$-sender for $\cT$ with signal edges $e_q$ and $e_{q+1}$. Then we are done using the previous case.\end{proof}

\begin{rem}
Using similar arguments, we can build $X$-senders and $X$-determiners for any set $X$ of the form $X=\bigcup_{i\in I}X_{i}$ with $X\subsetneq [q]$. These gadgets graphs are not required for our purposes but could be useful for other applications.
\end{rem}

\subsection{The basic construction}\label{sec:basic_construction}

\subsubsection{Special hypergraph}\label{sec:definition_hypergraph}
The construction described in this section is due to \Rodl{} and \Siggers{}~\cite[Section 7]{rodl_ramsey_2008}.

\begin{definition} Let $q,\ell\geq 2$ be integers.
\begin{enumerate}
    \item An \emph{oriented} $\ell$-graph $\mathbf{H}$ on vertex set $V=V(\mathbf{H})$ is a set of ordered $\ell$-element subsets of $V$, called \emph{arcs}.
    \item A \emph{$(q,\ell)$-pattern} is an $\ell$-tuple of not necessarily distinct elements from $[q]$ (i.e., an element of $[q]^{\ell}$).
    \item Given a set $\Phi\subsetneq [q]^{\ell}$ of forbidden $(q,\ell)$-patterns, a $q$-coloring $\chi:V(\mathbf{H})\to[q]$ of the vertices of $\mathbf{H}$ is \emph{$\Phi$-avoiding} if $\parens*{\chi(v_1),\ldots,\chi(v_\ell)}\not\in\Phi\text{ for any arc }(v_1,\ldots,v_\ell)\in E(\mathbf{H}).$
\end{enumerate}
\end{definition}

The following lemma establishes the existence of an oriented hypergraph $\mathbf{H}$
such that, in every $q$-coloring of the vertices of $\mathbf{H}$ that avoids a set of forbidden patterns, two special vertices are guaranteed to have distinct colors; the proof  can be found in \cite[Lemma 7.5]{rodl_ramsey_2008}.

\begin{lemma}[\Rodl{} and Siggers~\cite{rodl_ramsey_2008}]\thlabel{lemma:exist_hypergraph}
Given $q,\ell,g \geq 2$, let $\Phi\subsetneq [q]^{\ell}$ and assume $\Phi$ contains the monochromatic pattern $(i,i,\ldots,i)$ for each $i\in[q]$. Then there exists an oriented $\ell$-graph $\mathbf{H}=\mathbf{H}(\Phi,q,\ell,g)$ with girth more than $g$ and distinguished vertices $u$ and $u'$ satisfying the following properties:
\begin{enumerate}
    \item $\mathbf{H}$ has a $\Phi$-avoiding $q$-vertex coloring.
    \item There is no arc of $\mathbf{H}$ containing both vertices $u$ and $u'$. \label{lemma:hypergraph_RS_noarc}
    \item Under any $\Phi$-avoiding $q$-vertex coloring of $\mathbf{H}$, the vertices $u$ and $u'$ receive different colors.\label{lemma:hypergraph_different_colors}
\end{enumerate}
\end{lemma}

\subsubsection{The core construction}\label{sec:core_construction}

Here we describe the basic construction as given in~\cite{rodl_ramsey_2008}. However, as we will see later, for asymmetric tuples this construction does not necessarily produce a  useful gadget immediately. We will address this issue in the next subsection.
\smallskip

Let $q\geq 2$ and $\cT=(K_{t_1},\ldots,K_{t_q})$ be any $q$-tuple of cliques with $t_1\geq\ldots\geq t_q\geq 3$. Let $F'$ be a $q$-Ramsey-minimal graph for $\cT$ and $xy\in E(F')$, and let $F$ be obtained from $F'$ by removing the edge $xy$ (but keeping the endpoints). 

First note that, by the minimality of $F'$, there exists a $\cT$-free coloring of $F$.
We now claim that in any such coloring of $F$ the vertex $x$ must be incident to edges of every color. 
Indeed, let $\chi$ be a $\cT$-free coloring of $F$. For any $t\geq 3$, we have $\delta(K_t)\geq 2$, so if $x$ has no incident edge in some color $i$, we can extend $\chi$ to a coloring of $F'$ by setting $\chi(xy)=i$. Then $xy$ cannot participate in a new monochromatic clique, contradicting the fact that $F'$ is $q$-Ramsey for $\cT$. Therefore, the claim holds.

Let $\ell=d_F(x)$, and fix an orientation $(x_1,\ldots,x_\ell)$ of the neighborhood of $x$ in $F$. Let $\Phi$ be the set of $(q,\ell)$-patterns $(c_1,\ldots,c_\ell)$ such that the following holds: 
\begin{quote}
    The $q$-coloring $\chi$ of the edges $xx_1,\ldots,xx_\ell$ of $F$ defined by $\chi(xx_i)=c_i$ cannot be extended to a $\cT$-free coloring of $F$.
\end{quote}

Note that $\Phi$ contains all $(q,\ell)$-patterns with fewer than $q$ colors. Moreover, since $F$ is not $q$-Ramsey for $\cT$, there exists a $(q,\ell)$-pattern not in $\Phi$. Let $\bm{\mathcal{H}}$ be an oriented $\ell$-graph with girth at least four as guaranteed by \thref{lemma:exist_hypergraph} for the set $\Phi$. 

\begin{figure}[ht]
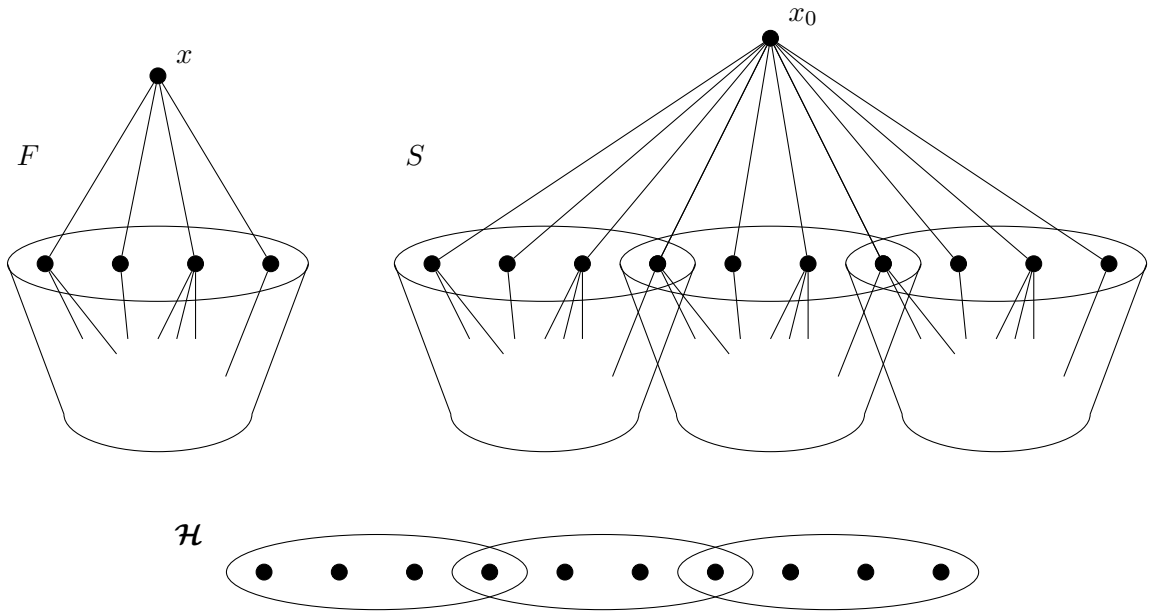

    \centering
    \tikzRS{0.5}\hspace{1cm}\vspace{0.8cm}
    \tikzRSFull{0.5}
    \tikzRSHyper{0.5}
    \caption{Core construction $F$, $S$, and the hypergraph $\mathbf{H}$, replicated from~\cite{rodl_ramsey_2008}.}
    \label{fig:core_construction}
\end{figure}

Let $S=S(\cT)$ be the following graph (see Figure~\ref{fig:core_construction} for an illustration). For any arc $\bm{a}\in\mathbf{H}$, let $F^{\bm{a}}$ be a copy of $F$ and $x^{\bm{a}}$ be the copy of $x$ in $F^{\bm{a}}$ and identify the vertices of $\bm{a}$ with those of $N(x^{\bm{a}})$, respecting the previously fixed orientation of $N_F(x)$. Finally, identify the vertices $x^{\bm{a}}$ for all arcs $\bm{a}$ and call the resulting vertex $x_0$.

\subsubsection{Properties of the basic construction}\label{sec:core_properties}

We write $e= x_0u$ and $f=x_0u'$.
\begin{claim}\thlabel{claim:propertiesS_new}
The graph $S$ has the following properties:
\begin{enumerate}[label=(\alph*)]
    \item $S$ is not $q$-Ramsey for $\cT$. \label{claim:S_not_ramsey_new}
    \item In any $\cT$-free coloring $\phi$ of $S$, the edges $e$ and $f$ receive different colors. \label{claim:S_different_colors_new}
    \item The edges $e$ and $f$ are adjacent. \label{claim:S_ef_adjacent_new}
    \item There is no triangle containing $e$ and $f$ in $S$.\label{claim:S_no_triangle_new}
\end{enumerate}
\end{claim}

The proofs of~\ref{claim:S_not_ramsey_new} and~\ref{claim:S_different_colors_new} are given in Claims 7.7, 7.8, and 7.9 in~\cite{rodl_ramsey_2008}. Part~\ref{claim:S_no_triangle_new} follows directly from \thref{lemma:exist_hypergraph}\ref{lemma:hypergraph_RS_noarc}: As there is no arc of $\mathbf{H}$ containing $u$ and $u'$, there is no edge between $u$ and $u'$ in $S$.
\medskip

Now, if the tuple $\cT$ is symmetric, then this construction yields a negative $[q]$-sender for $\cT$. Indeed, Property~\ref{axiom:set_sender_not_Ramsey} is guaranteed by~\thref{claim:propertiesS_new}\ref{claim:S_not_ramsey_new}, Property~\ref{axiom:set_sender_colored_edges} is given by~\thref{claim:propertiesS_new}\ref{claim:S_different_colors_new}, and Property~\ref{axiom:set_sender_any_color} is true also by~\thref{claim:propertiesS_new}\ref{claim:S_not_ramsey_new} since we can always permute colors. This negative $[q]$-sender can then be used to construct safe positive and negative $[q]$-senders using standard constructions (see~\cite{rodl_ramsey_2008}). Unfortunately, in the asymmetric setting~\thref{claim:propertiesS_new}\ref{claim:S_not_ramsey_new} does not allow us to conclude Property~\ref{axiom:set_sender_any_color}. If there exists a color $c\in [q]$ such that $x_0u$ does not have color $c$ in any $\cT$-free coloring, then $S$ is a safe $X$-determiner with signal edge $x_0u$ for some $\emptyset \subsetneq X\subseteq [q]\setminus\set{c}$. A similar argument applies with $x_0u'$ instead of $x_0u$. In either case, we are done by~\thref{thm:one_det_is_enough}. But it is possible that each of the $q$ colors can appear on each of $x_0u$ and $x_0u'$ in a $\cT$-free coloring but Property~\ref{axiom:set_sender_any_color} still fails to hold. Thus, we will need to take the construction a step further to overcome this challenge. This is the content of the next subsection.

\subsection{Building set-senders and set-determiners}\label{sec:constructions_new}

We begin with some general statements and a definition that will be useful for our subsequent arguments. First we prove the following useful technical lemma, which allows us to put together two graphs and color each independently with a $\cT$-free coloring to obtain a $\cT$-free coloring of the larger graph (cf.\ safeness of set-senders and determiners). 

\begin{figure}[ht!]
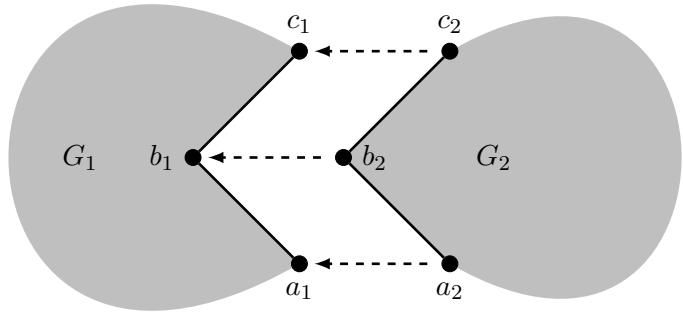

    \centering
    \tikzjoinS{2}
    \caption{Identifying the edges of two graphs.}
    \label{fig:identifying}
\end{figure}

\begin{lemma}\thlabel{lem:safeness_new}
Let $\cT$ be any $q$-tuple of cliques. Let $G_1$ and $G_2$ be two graphs with disjoint sets of vertices. For each $i\in [2]$, let $a_i,b_i,c_i\in V(G_i)$ be vertices such that $a_ib_i, b_ic_i\in E(G_i)$ but $a_ic_i\not\in E(G_i)$.
Let $G$ be the (simple) graph obtained by identifying the following pairs of vertices: $a_1$ with $a_2$, $b_1$ with $b_2$, and $c_1$ with $c_2$ (see Figure~\ref{fig:identifying}). Call the newly obtained vertices $a, b$, and $c$, respectively. Then a $q$-coloring $\phi$ of $G$ is $\cT$-free if and only if $\phi$ induces a $\cT$-free coloring on each of $G_1$ and $G_2$. In particular, $G$ has a $\cT$-free coloring if and only if $G_1$ and $G_2$ have $\cT$-free colorings $\phi_1$ and $\phi_2$ satisfying $\phi_1(a_1b_1) = \phi_2(a_2b_2)$ and $\phi_1(b_1c_1) = \phi_2(b_2c_2)$.
\end{lemma}

\begin{proof}
We show the first part of the statement. The ``in particular'' part then follows immediately.
If $\phi$ is a $\cT$-free coloring of $G$, then it clearly induces $\cT$-free colorings on $G_1$ and $G_2$. We now prove the other direction. Let $\phi$ be a coloring of $G$ and assume $\phi_{\mid G_1}$ and $\phi_{\mid G_2}$ are both $\cT$-free. Suppose for a contradiction that $\phi$ is not $\cT$-free on $G$, that is, suppose there exists a color $i\in [q]$ and a copy $K$ of $K_{t_i}$ in color $i$ in $G$. Then $K$ must contain at least one edge from $E(G_1)\setminus E(G_2)$ and one edge from $E(G_2)\setminus E(G_1)$. Since by assumption neither $G_1$ nor $G_2$ contains an edge between the vertices $a$ and $c$, the clique $K$ must contain at least one vertex from each of $V(G_1)\setminus V(G_2)$ and $V(G_2)\setminus V(G_1)$. But two such vertices cannot share an edge in $G$ and thus cannot be part of $K$, a contradiction.
\end{proof}

For a graph $G$ with distinguished edges $g_1$ and $g_2$,  
we will use a digraph to encode which pairs of colors can occur on the edges $g_1$ and $g_2$ in a $\cT$-free coloring of $G$.

\begin{definition}[Auxiliary digraph $D(\cT,G,g_1,g_2)$]\thlabel{def:D_new}
Let $\cT$ be a $q$-tuple of cliques and $G$ be a graph with two distinguished edges $g_1$ and $g_2$. We define $D=D(\cT,G,g_1,g_2)$ to be a digraph on vertex set $[q]$ in which, for any $i,j\in[q]$, the arc $\vv{ij}$ is in $D$ if and only if there exists a $\cT$-free coloring of $G$ in which $g_1$ receives color $i$ and $g_2$ receives color $j$.
\end{definition}

For example, if $G$ is Ramsey for $\cT$ then, for any edges $e$ and $f$, the digraph $D=D(\cT,G,e,f)$ contains no arcs. If $G$ is a negative $X$-sender for $\cT$ with signal edges $e$ and $f$, then $D=D(\cT,G,e,f)$ consists of a complete directed subgraph on $X$ (with arcs in both directions between any pair of vertices in $X$ but without loops), and $q-\size{X}$ isolated vertices. If $G$ is a positive $X$-sender for $\cT$ with signal edges $e$ and $f$, then $D=D(\cT,G,e,f)$ contains a loop at each vertex $i\in X$ and no other arcs.\medskip

In the remainder of this section, to avoid repetitions, $\cT$ always denotes a $q$-tuple of cliques $\cT=\parens*{K_{t_1},\ldots,K_{t_q}}$ with $q\geq 3$ and $t_1\geq\ldots\geq t_q\geq 3$, $S=S(\cT)$ denotes the gadget $S$ as built in Section~\ref{sec:basic_construction}, with distinguished edges $e = x_0u$ and $f = x_0u'$, and $D=D(\cT,S,e,f)$ denotes the auxiliary digraph of $S$ as constructed in \thref{def:D_new}. This definition implies that:
\begin{enumerate}
    \item In any $\cT$-free coloring $\phi$ of $S$,  if $\phi(e)=i$, then $\phi(f)\in N^+_D(i)$.
    \item For any colors $i,j$ such that $j\in N^+_D(i)$, there exists a $\cT$-free coloring $\phi$ of $S$ such that $\phi(e)=i$ and $\phi(f)=j$.\label{claim:S_any_colors_new}
\end{enumerate}\medskip

As explained in Section~\ref{sec:basic_construction}, if there exists a color that cannot appear on one of the distinguished edges of $S$ in any $\cT$-free coloring, then $S$ is already a set-determiner. The next lemma is a reformulation of this idea in the language of the auxiliary digraph $D$.

\begin{lemma}[Basic conditions]\thlabel{lem:basic_conditions_new}
Suppose that there exists a vertex $i\in [q]$ of $D$ such that the in- or out-neighborhood of $i$ in $D$ is empty,  that is, $N^-_{D}(i)=\emptyset$ or $N^+_{D}(i)=\emptyset$. Then there exists a safe $X$-determiner for $\cT$  for some non-empty $X\subsetneq [q]$.
\end{lemma}
\begin{proof}
    Suppose the out-neighborhood of $i$ is empty. This means that $e$ cannot receive color $i$ in any $\cT$-free coloring of $S$. As $S$ is not $q$-Ramsey for $\cT$, it follows that $S$ is an $X$-determiner with signal edge $e$, for some non-empty subset of colors $X\subseteq[q]\setminus\set{i}$. By \thref{thm:all_det_are_safe}, it is furthermore a safe determiner. A similar argument applies when the in-neighborhood of $i$ in $D$ is empty with $f$ instead of $e$.
\end{proof}

First we show that if $D$ contains no $2$-cycle, we can construct a suitable determiner.

\begin{lemma}[Existence of a $2$-cycle]\thlabel{lem:existence_2_cycle}
If $D$ contains no $2$-cycle, then there exists a safe $X$-determiner for $\cT$ for some $X\subsetneq [q]$.
\end{lemma}

\begin{proof}
Suppose that $D$ contains no $2$-cycles, in particular there is no $j\in[q]$ such that $j\in N^+_{D}(1)$ and $j\in N^-_{D}(1)$.

Let $A$ be the second out-neighborhood of $1$ in $D$, that is, let $A$ be given by
\[A=\set*{u\in [q]:\ \exists v\in [q],\ \vv{1v}\in E(D) \text{ and }\vv{vu}\in E(D)}.\] We first claim that $A$ contains no color corresponding to a clique of size $t_1$, that is, $\set*{j : t_j=t_1}\cap A = \emptyset$. Suppose there exists  $j\in A$ with $t_j=t_1$. Thus, there exists a directed $2$-path from $1$ to $j$ in $D$. Let $k$ denote the intermediate vertex in one such path. We know that there exists a $\cT$-coloring in which $e$ receives color $k$ and $f$ receives color $j$. As $t_j=t_1$, by switching the roles of colors $1$ and $j$ in this coloring, we obtain another $\cT$-free coloring of $S$ implying that $\vv{k1}$ is also an arc of $D$, creating a $2$-cycle $\vv{1k1}$ in $D$, a contradiction.

Since the largest clique in $\cT_{A}$ has size at most $t_1-1$, using \thref{thm:NR}, there exists a $K_{t_1}$-free graph $H$ such that $H\to \cT_{A}$. Without loss of generality, let $V(H)=[m]$, and let $u$ and $v$ be distinct vertices disjoint from $H$. Add an edge between $v$ and $u$, and for every $x\in [m-1]$ add an edge between $v$ and $x$. For every edge $h=xy$ in $E(H)$, with $1\leq x<y\leq m$, let $S_1,S_2$ be two copies of $S$ with distinguished edges $e_1,f_1$ and $e_2,f_2$ respectively, and 
\begin{itemize}
    \item join $uv$ and $vx$ by $S_1$, identifying $uv$ with $e_1$ and $vx$ with $f_1$,
    \item join $vx$ and $h$ by $S_2$, identifying $vx$ with $e_2$ and $h$ with $f_2$.
\end{itemize}
Outside of the specified intersections, all copies of $S$ are vertex- and edge-disjoint. Call the resulting graph $G$. For clarity, for every $y>x$, there exists a copy $S_1$ of $S$ sitting on $uv$ and $vx$ and  a copy $S_2$ of $S$ sitting on $vx$ and $xy$. Therefore, for every edge $h$ in $H$, the graph $G$ includes the subgraph displayed in Figure~\ref{fig:graphS2_new}. 

\begin{figure}[ht!]
    \centering
    \tikzSTwoNew{2}
    \caption{Proof of \thref{lem:existence_2_cycle}, subgraph of $G$ at every edge $h$ in $H$.}
    \label{fig:graphS2_new}
\end{figure}

We claim that $G$ is an $X$-determiner with signal edge $uv$ for some $X\subsetneq [q]$. By \thref{lem:basic_conditions_new}, we may assume that there exist $k\in N^-_D(1)$ and $j\in N^-_D(k)$. By~\thref{claim:propertiesS_new}\ref{claim:S_different_colors_new}, $D$ contains no loops, therefore we know that $k\neq 1$ and $k\neq j$, and as $D$ contains no $2$-cycle we know that $1\neq j$. Therefore $D$ contains the directed path $\vv{jk1}$. In particular, we will show $1\not\in X$ and $j\in X$.

First we show that $G\nrightarrow_q \cT$. Let $\phi$ be the following coloring of $G$. Set $\phi(uv) = j$. For every edge $vx$ with $x\in V(H)$, set $\phi(vx)=k$, and for any edge $h\in E(H)$, set $\phi(h)=1$.
As $k\in N^+_D(j)$ and $1\in N^+_D(k)$, we can extend this coloring to the copies of $S$ so that each of these copies receives a $\cT$-free coloring. We claim that this yields a $\cT$-free coloring of all of $G$. Suppose not; then there exists a color $\ell\in [q]$ and a copy $K$ of $K_{t_\ell}$ that is monochromatic in color $\ell$. Now, if $K$ contains a vertex not in $V(H)\cup \set{u,v}$, this vertex must be contained within some copy of $S$ and because $K$ is a clique, all of its neighbors must be within this particular copy of $S$. This in turn implies that $K$ cannot be monochromatic because each copy of $S$ was given a $\cT$-free coloring, a contradiction. 

Therefore $V(K)\subseteq V(H)\cup\set*{u,v}$. As $H$ is monochromatic in color $1$ and $K_{t_1}$-free and we added no further edges between the vertices of $H$, it follows that $V(K)\cap\set*{u,v}\neq\emptyset$. Note that by construction the vertex $u$ does not have any neighbors in $V(H)$, so it cannot be part of $K$. 
If $v\in V(K)$, as any edge from $v$ to $V(H)$ corresponds to an edge $vx$ as shown in Figure~\ref{fig:graphS2_new}, then $K$ must use at least one such edge and hence must be monochromatic in color $k\neq 1$.
As $H$ is monochromatic in color $1$, $K$ cannot contain any edge from $H$, so it must be fully contained in the star between $v$ and $V(H)$. This is, however, not possible since $t_\ell\geq 3$. Therefore $\phi$ is a $\cT$-free coloring of $G$ with $\phi(uv)=j$.

Now, assume that there exists a $\cT$-free coloring of $G$ with $\phi(uv)=1$. By~\thref{def:D_new}, every edge of $H$ is colored with a color in $A$, the second out-neighborhood of $1$ in $D$. As $H\to \cT_{A}$, this coloring is not $\cT$-free, a contradiction.

Therefore $G$ is an $X$-determiner for $\cT$ for some $X$ with $j\in X$ and $1\not\in X$, as required. By \thref{thm:all_det_are_safe}, it is furthermore a safe determiner. \end{proof}

The next construction allows us to build a \emph{symmetric} gadget, that is, a gadget for which each vertex has identical in- and out-neighborhood in its auxiliary digraph. This will considerably simplify all future arguments.

\begin{figure}[ht!]
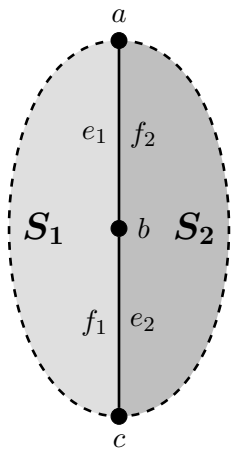

    \centering
    \tikzSymmetricSVariant{2.5}
    \caption{The graph $S^\sym$.}
    \label{fig:graphS_Sinverse_var}
\end{figure}

\begin{definition}[Symmetric construction $S^\sym$]\thlabel{def:S_symmetric_new}
Let $S_1$ and $S_2$ be two copies of $S$, with respective distinguished edges $e_1,f_1$ and $e_2,f_2$. We define $S^\sym=S^\sym(\cT)$ to be the following graph. Let $abc$ be a path on three vertices. Join $ab$ and $bc$ by $S_1$, identifying $ab$ with $e_1$ and $bc$ with $f_1$. Join $ab$ and $bc$ by $S_2$, identifying $ab$ with $f_2$ and $bc$ with $e_2$. See Figure~\ref{fig:graphS_Sinverse_var} for an illustration. We let $ab$ and $bc$  be the distinguished edges of $S^\sym$, and we further define $D^\sym=D(\cT,S^\sym,ab,bc)$ to be the auxiliary graph of $S^\sym$.
\end{definition}

The next lemma shows that if $D$ contains at least one $2$-cycle, then $S^\sym$ satisfies all properties in~\thref{claim:propertiesS_new}, while also being symmetric.

\begin{lemma}\thlabel{lem:Ssym_properties_new}
Let $S^\sym$ and $D^\sym$ be as in~\thref{def:S_symmetric_new}, with $ab$ and $bc$ being the distinguished edges of $S^\sym$. If $D$ contains at least one $2$-cycle, then $S^\sym$ satisfies the following properties:
    \begin{enumerate}[label=(\alph*)]
        \item $S^\sym$ is not $q$-Ramsey for $\cT$. \label{claim:Ssym_not_ramsey_new}
        
        \item In any $\cT$-free coloring of $S^\sym$, the edges $ab$ and $bc$ are assigned different colors. \label{claim:Ssym_different_colors_new}
        
        \item The edges $ab$ and $bc$ are adjacent. \label{claim:Ssym_ef_adjacent_new}
        
        \item There is no triangle containing $ab$ and $bc$ in $S^\sym$.\label{claim:Ssym_no_triangle_new}
        
        \item For any pair of colors $i,j\in [q]$, there exists a $\cT$-free coloring $\phi$ of $S^\sym$ with $\phi(ab)=i$ and $\phi(bc)=j$ if and only if there exists a $\cT$-free coloring $\psi$ of $S^\sym$ with $\psi(ab)=j$ and $\psi(bc)=i$. That is, $\vv{ij}\in E(D^\sym)$ if and only if $\vv{ji}\in E(D^\sym)$.\label{claim:Ssym_symmetric_new}
    \end{enumerate}
\end{lemma}

\begin{proof}
For Property~\ref{claim:Ssym_not_ramsey_new}, as $D$ contains at least one $2$-cycle, let $i,j\in [q]$ be such that $\vv{ij}\in E(D)$ and $\vv{ji}\in E(D)$. Then there exists a $\cT$-free coloring of $S_1\cong S$ in which $e_1$ receives color $i$ and $f_1$ receives color $j$ and a $\cT$-free coloring of $S_2\cong S$ in which $e_2$ receives color $j$ and $f_2$ receives color $i$. These two colorings agree on the edges $ab$ and $bc$ of $S^\sym$, so~\thref{lem:safeness_new} allows us to conclude that there exists a $\cT$-free coloring of $S^\sym$. 

Property~\ref{claim:Ssym_no_triangle_new} follows easily from~\thref{claim:propertiesS_new}\ref{claim:S_no_triangle_new}. Property~\ref{claim:Ssym_ef_adjacent_new} is trivial by construction. Property~\ref{claim:Ssym_symmetric_new} follows from the following fact: by~\thref{lem:safeness_new}, $S^\sym$ has a $\cT$-free coloring in which $ab$ has color $i$ and $bc$ has color $j$ if and only if $S_1$ has a $\cT$-free coloring in which $e_1$ has color $i$ and $f_1$ has color $j$ and $S_2$ has a $\cT$-free coloring in which $e_2$ has color $j$ and $f_2$ has color $i$; but $S_1$ and $S_2$ are copies of the same graph, so we can exchange these two colorings to obtain the required second coloring of $S^\sym$.

Finally Property~\ref{claim:Ssym_different_colors_new} is also not difficult to see: any $\cT$-free coloring of $S^\sym$ induces a $\cT$-free coloring on the edges of $S_1$ and by~\thref{claim:propertiesS_new}\ref{claim:S_different_colors_new}, we know that the edges $e_1=ab$ and $f_1=bc$ must receive different colors.\end{proof}

\begin{rem}\thlabel{rem:symmetry}
Observe that Properties~\ref{claim:Ssym_not_ramsey_new}--\ref{claim:Ssym_no_triangle_new} of this ``symmetrized'' gadget $S^\sym$ are exactly all the properties of $S$ from~\thref{claim:propertiesS_new}. Therefore, whenever we are able to assume that $D$ contains at least one $2$-cycle, we will (slightly abusing notation) assume that $S$ is itself symmetric, that is,~\thref{lem:Ssym_properties_new}\ref{claim:Ssym_symmetric_new} is true for $S$. This in particular means that, for any vertex $i\in V(D)$, we have $N^+_D(i)=N^-_D(i)$. When this holds, for convenience we will write $N_D(i)$ instead of $N^+_D(i)$ and $N^-_D(i)$.
\end{rem}

The following construction will be useful in the subsequent proof. 

\begin{definition}[Claw construction]\thlabel{def:claw_construction_new}
For given integers $h,d$ such that $h>d$, let $\claw_\cT(h,d)$ be the following graph. Let $H$ be a copy of $K_h$ on vertex set $\set*{v_1,\ldots,v_h}$ and $xy$ be a disjoint edge. Add edges from $x$ to the vertices $v_1,\ldots,v_{d}$ of $H$. For each $i\in[d]$, join the edges $xy$ and $xv_i$ by a copy of $S$, identifying $xy$ with $e$ in $S$ and $xv_i$ with $f$ in $S$. Finally add an edge from $x$ to $v_{d+1}$. See Figure~\ref{fig:claw_construction_new} for an illustration.
\end{definition}

\begin{figure}[ht]
    \centering
    \tikzclawComplement{0.8}
    \caption{The Claw construction $\claw_\cT(7,3)$.}
    \label{fig:claw_construction_new}
\end{figure}

The Claw construction allows us to build a gadget graph, very similar to $S$, but with ``opposite'' pairs of colors allowed on its distinguished edges.

\begin{lemma}\thlabel{lem:Sc_new}
    At least one of the following statements is true:
    \begin{enumerate}
        \item There exists an $X$-determiner for $\cT$ for some non-empty $X\subsetneq [q]$.\label{lem:Sc_det_new}
        \item There exist parameters $h,d$ such that $S^c=\claw_\cT(h,d)$ is a gadget with distinguished edges $xy$ and $xz$ satisfying the following properties:\label{lem:Sc_existence_new}
        \begin{enumerate}
            \item $S^c$ is not $q$-Ramsey for $\cT$. \label{lem:Sc_not_ramsey_new}
            \item For every $\cT$-free coloring of $S^c$ and for every $i\in[q]$, if $xy$ receives color $i$, then $xz$ receives a color in $[q]\setminus N^+_D(i)$.\label{lem:Sc_colored_edges_new}
            \item For any colors $i,j$ such that $j\in [q]\setminus N^+_D(i)$, there exists a $\cT$-free coloring of $S^c$ such that $xy$ is colored $i$ and $xz$ is colored $j$.\label{lem:Sc_any_color_new}

        \end{enumerate}
    \end{enumerate}
\end{lemma}
\begin{proof}

    First note that, if $D$ does not contain any $2$-cycle, then we are done by \thref{lem:existence_2_cycle}. Therefore, we can assume that $D$ contains at least one $2$-cycle and hence, by \thref{lem:Ssym_properties_new,rem:symmetry}, that $S$ is symmetric, that is, \thref{lem:Ssym_properties_new}\ref{claim:Ssym_symmetric_new} is true for $S$.
    
    Let $h=r_q(\cT)-1$, and for each $0\leq d\leq h-1$, let $G_d=\claw_{\cT}(h,d)$. Note that, in any $\cT$-free coloring of $G_d$ (if such a coloring exists), if the edge $xy$ has color $i\in[q]$, then each edge $xv_j$ for $0\leq j \leq d$ in $G_d$ must be colored with a color in $N_D(i)$. 
    
    It is not difficult to check that $G_0$ is not Ramsey for $\cT$: indeed, the choice of $h$ implies that $K_h$ has a $\cT$-free $q$-coloring, and $G_0$ consists of $K_h$ and a pendent path $yxv_1$ of length two. Since every clique in $\cT$ has size at least $3$, this pendent path, irrespective of the coloring of its edges, cannot create a monochromatic clique. On the other hand, $G_{h-1}$ is Ramsey for $\cT$, as it contains a clique on $r_q(\cT)$ vertices. 
    
    Let $i\in [q]$ be fixed. By the previous paragraph, we know that $G_0$ has a $\cT$-free coloring in which $xy$ receives color $i$, while $G_{h-1}$ has no such coloring. Let $0\leq k_i<h-1$ be the maximum integer such that there exists a $\cT$-free coloring of $G_{k_i}$ in which the edge $xy$ is colored $i$. Now, if $k_i<k_j$ for some $i,j\in [q]$, then $G_{k_j}$ is an $X$-determiner for $\cT$ with signal edge $xy$, for some $X$ with $j\in X$ and $i\not\in X$, hence verifying~\ref{lem:Sc_det_new}. \medskip
    
    Assume now that $k_1=\dots = k_q = k$ for some integer $0\leq k<h-1$, and let $z=v_{k+1}$. Set $S^c = G_{k}$. We claim that this graph satisfies the required properties of~\ref{lem:Sc_existence_new}. 
    
    By the definition of $k$, for any $i\in[q]$, we have $k_i=k$ and there exists a $\cT$-free coloring of $S^c$ in which the edge $xy$ is colored $i$, hence verifying Property~\ref{lem:Sc_not_ramsey_new}. 
    
    Assume that there exists a $\cT$-free coloring $\phi$ of $S^c$ such that $xy$ receives color $i\in[q]$ and $xv_{k+1}$ receives color $j\in N_D(i)$. Note that $G_{k+1}$ can be obtained from $G_k$ by joining the edges $xy$ and $xv_{k+1}$ by a copy of $S$ and adding the edge $xv_{k+2}$ (see Figure~\ref{fig:claw_construction_induction}). 
    Apply $\phi$ as a partial-coloring of $G_{k+1}$. By the definition of $i,j$, and by \thref{lem:safeness_new}, we can extend $\phi$ to the copy of $S$ joining the edges $xy$ and $xv_{k+1}$, forming a $\cT$-free coloring of $G_{k+1}-{xv_{k+2}}$. 
    
    Set $\phi(xv_{k+2})=i$. Recall that by~\thref{claim:propertiesS_new}\ref{claim:S_different_colors_new}, we know that $i\not\in N_D(i)$.
    By construction, the common neighbors of $x$ and $v_{k+2}$ are all in $V(H)$. However, each edge $xv_\ell$ for $0\leq \ell \leq k+1$ must be colored with a color in $N_D(i)$, while $i\not\in N_D(i)$.
    As $t_i\geq 3$, the new edge $xv_{k+2}$ cannot be part of a monochromatic clique $K_{t_i}$ in color $i$, and therefore the coloring $\phi$ is a $\cT$-free coloring of $G_{k+1}$, contradicting the maximality of $k$. Therefore in any $\cT$-free coloring $\phi$ of $S^c$ such that $xy$ receives color $i\in[q]$, the edge $xv_{k+1}$ receives some color $j\in [q]\setminus N_D(i)$, verifying Property~\ref{lem:Sc_colored_edges_new}.

\begin{figure}[ht]
    \centering
    \tikzclawComplementInduction{0.8}
    \caption{Proof of \thref{lem:Sc_new}, from $G_2$ to $G_3$.}
    \label{fig:claw_construction_induction}
\end{figure}

    Finally, to prove~\ref{lem:Sc_any_color_new}, let $i,j\in [q]$ be such that $j\not\in N_{D}(i)$. Let $\phi$ be a $\cT$-free partial coloring of $S^c$ such that $xz$ is the only uncolored edge and $xy$ receives color $i\in[q]$, and set $\phi(xz)=j$. 
    Similarly to the above argument, the only common neighbors of $x$ and $z$ are vertices $v_\ell$ from $H$ with $0\leq \ell \leq k$. However, each edge $xv_\ell$ for $0\leq \ell \leq k$ must be colored with a color in $N_D(i)$, but $j\not\in N_D(i)$.
    Therefore, as $t_j\geq 3$, this is still a $\cT$-free coloring of $S^c$.\end{proof}

From now on, whenever we refer to the gadget $S^c$, we will assume that~\thref{lem:Sc_new}\ref{lem:Sc_existence_new} holds. We will further write $D^c = D(\cT, S^c, xy, xz)$ for its auxiliary digraph. Note that, for all $i\in [q]$, we have $N^+_D(i)=N^-_D(i)$ if and only if $N^+_{D^c}(i)=N^-_{D^c}(i)$. If these equalities hold, we will denote $N^+_{D^c}(i)$ and $N^-_{D^c}(i)$ by $N_{D^c}(i)$. The following observation will be useful in the subsequent proofs.

\begin{obs}\thlabel{obs:loop_each_vertex}
    For any $i\in[q]$, we have $N^+_{D^c}(i)=[q]\setminus N^+_D(i)$. In particular $D^c$ contains a loop at every vertex and therefore $i\in N^+_{D^c}(i)\neq\varnothing$ for all $i\in[q]$.
\end{obs}

We now present the final two intermediate constructions, obtained by joining copies of $S$ and $S^c$ in a sequential fashion.

\begin{figure}[ht]
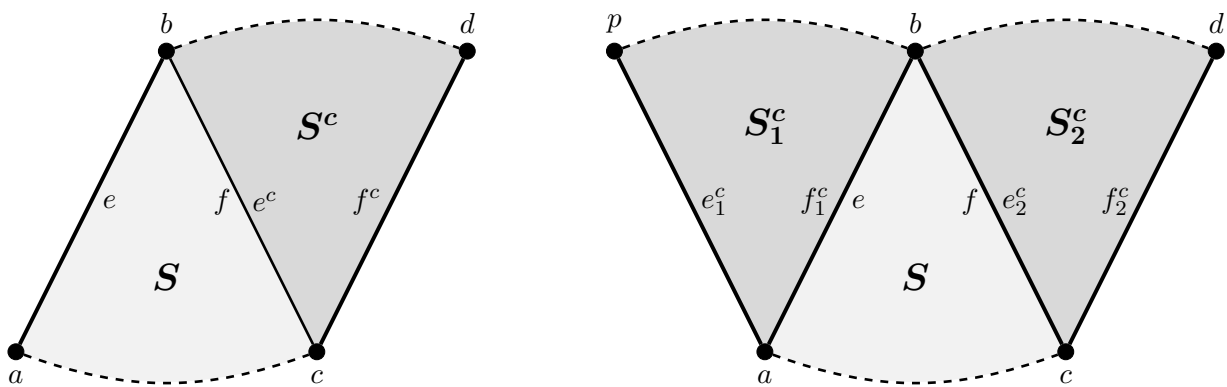

    \centering
    \tikzSScomplement{2}\hspace{1.25cm}
    \tikzScSSc{2}
    \caption{Constructions $T$ (left) and $T'$ (right).}
    \label{fig:S_Scomplement_new}
\end{figure}

\begin{definition}[Construction $T$]\thlabel{def:T_new}
Let $S^c$ be the gadget from~\thref{lem:Sc_new} with distinguished edges $e^c=xy$ and $f^c=xz$. Let $T$ denote the graph obtained from $S$ and $S^c$ as follows. Let $abcd$ be a path with three edges. Join $ab$ and $bc$ by a copy of $S$, identifying $ab$ and $bc$ with the copies of $e$ and $f$ (respectively). Join $bc$ and $cd$ by a copy of $S^c$, identifying $bc$ and $cd$ with the copies of $e^c$ and $f^c$ (respectively). See Figure~\ref{fig:S_Scomplement_new} (left) for an illustration. We let $ab$, $bc$ and $cd$ be the distinguished edges of $T$, and we further define $B=D(\cT,T,ab,cd)$ to be the auxiliary digraph of $T$.
\end{definition}

\begin{definition}[Construction $T'$]\thlabel{def:Tprime_new}
Let $S^c$ be the gadget from~\thref{lem:Sc_new} with distinguished edges $e^c=xy$ and $f^c=xz$. Let $T'$ denote the graph obtained from $S$ and $S^c$ as follows. Let $pabcd$ be a path with four edges. Join $pa$ and $ab$ by a copy $S^c_1$ of $S^c$, identifying $pa$ and $ab$ with the copies of $e^c$ and $f^c$ (respectively). Join $ab$ and $bc$ by a copy of $S$, identifying $ab$ and $bc$ with the copies of $e$ and $f$ (respectively). Join $bc$ and $cd$ by a copy $S^c_2$ of $S^c$, identifying $bc$ and $cd$ with the copies of $e^c$ and $f^c$ (respectively). See Figure~\ref{fig:S_Scomplement_new} (right) for an illustration. We let $pa$, $ab$, $bc$, and $cd$ be the distinguished edges of $T'$, and we further define $B'=D(\cT,T',pa,cd)$ to be the auxiliary digraph of~$T'$.
\end{definition}

Note that, although one could show it to be true, there is no immediate reason to assume that $T$ and $T'$ are symmetric, and hence that  the in- and out-neighborhoods of all vertices in each of $B$ and $B'$ are identical. As this fact is not required in our proof, we omit its proof.

\begin{lemma}[Properties of $T$ and $T'$]\thlabel{lem:SSc_properties_new}
    There exists an $X$-determiner for $\cT$ for some non-empty $X\subsetneq [q]$, or the graphs $T$ and $T'$ as in \thref{def:T_new,def:Tprime_new} exist and satisfy the following properties:
    
    \begin{enumerate}
        \item $T$ and $T'$ are not $q$-Ramsey for $\cT$.\label{lem:SSc_not_ramsey_new}

        \item For every $\cT$-free coloring of $T$, $ab$ and $cd$ receive different colors.\label{lem:SSc_colored_edges_new}
    
        \item For any colors $i\in [q]$, $j\in N^+_D(i)$, and $k\in N^+_{D^c}(j)$, there exists a $\cT$-free coloring of $T$ such that $ab,bc,cd$ are colored $i,j,k$, respectively. \label{lem:SSc_any_color_new}
        
        \item For any colors $i\in [q]$, $j\in N^+_{D^c}(i)$, $k\in N^+_D(j)$, $\ell\in N^+_{D^c}(k)$, there exists a $\cT$-free coloring of $T'$ such that $pa,ab,bc,cd$ are colored $i,j,k,\ell$, respectively. \label{lem:ScSSc_any_color_new}    
        
    \end{enumerate}
\end{lemma}

\begin{proof}
    As before, note that, if $D$ does not contain any $2$-cycle, then we are done by \thref{lem:existence_2_cycle}. Therefore, we can assume that $D$ contains at least one $2$-cycle and hence, by \thref{lem:Ssym_properties_new,rem:symmetry}, that $S$ is symmetric. 
    
    We first show part~\ref{lem:SSc_any_color_new}. Let $i\in [q]$ be any color, $j\in N_D(i)$, and $k\in N_{D^c}(j)$. Then, by the definition of $D$ and by~\thref{lem:Sc_new}\ref{lem:Sc_any_color_new}, we can extend this coloring to the copies of $S$ and $S^c$ in $T$ so that each receives a $\cT$-free coloring. Using~\thref{lem:safeness_new} twice (first adding the copy of $S$ to the path and then that of $S^c$ to the resulting graph), this is a $\cT$-free coloring of $T$. Part~\ref{lem:ScSSc_any_color_new} follows similarly. For part~\ref{lem:SSc_not_ramsey_new}, we know using~\thref{lem:basic_conditions_new} that for any $i\in[q]$ there exists $j\in N^+_{D}(i)$, and by~\thref{obs:loop_each_vertex} we have $j\in N^+_{D^c}(j)$. Therefore, it follows directly from~\ref{lem:SSc_any_color_new} that there exists a $\cT$-free coloring of $T$, and similarly for $T'$ using~\ref{lem:ScSSc_any_color_new}.
    
    For part~\ref{lem:SSc_colored_edges_new}, consider a $\cT$-free coloring of $T$ and let $i$, $j$, and $k$ be the colors of the edges $ab$, $bc$, and $cd$, respectively. By the properties of $S$ and $S^c$, we must have $j\in N_D(i)$ and $k\in [q]\setminus N_D(j)$, or equivalently $i\in N_D(j)$ and $k\not\in N_D(j)$, and hence $i\neq k$.\end{proof}

\begin{lemma}\thlabel{lem:final_B_complete_new}
There exists a safe $X$-determiner for $\cT$ for some non-empty $X\subsetneq [q]$.
\end{lemma}

\begin{proof}

Once again, we can assume that $D$ contains at least one $2$-cycle and hence, by \thref{lem:Ssym_properties_new,rem:symmetry}, that $S$ is symmetric. Further, we can assume that~\thref{lem:Sc_new}\ref{lem:Sc_existence_new} holds, as otherwise~\thref{lem:Sc_new}\ref{lem:Sc_det_new} holds and we are done.

Let $B$ and $B'$ denote the auxiliary digraphs of $T$ and $T'$ (respectively) as introduced in \thref{def:T_new,def:Tprime_new}. By~\thref{lem:SSc_properties_new}\ref{lem:SSc_colored_edges_new} it follows that, for any $i\in[q]$, we have $i\not\in N^+_{B}(i)$. In particular $1\not\in N^+_{B}(1)$. As in the proof of \thref{lem:existence_2_cycle}, we first claim that $N^+_{B}(1)$ contains no color corresponding to a clique of size $t_1$. Suppose there exists some $j\in N^+_{B}(1)$ with $t_j=t_1$. Thus, there exists $k\in[q]$ such that $k\in N_D(1)$ and $j\in N_{D^c}(k)$. Since $k\in N_D(1)$, there exists a $\cT$-free coloring of $S$ in which $e$ receives color $1$ and $f$ receives color $k$. As $t_j=t_1$, by switching the roles of colors $1$ and $j$ in this coloring, we obtain another $\cT$-free coloring of $S$ in which $e$ receives color $j$ and $f$ receives color $k$; hence $k\in N_D(j)$, which in turn implies $j\in N_D(k)$ and $j\not\in N_{D^c}(k)$, a contradiction. Therefore the largest clique in $\cT_{N^+_{B}(1)}$ has size at most $t_1-1$, and by~\thref{thm:NR}, there exists a $K_{t_1}$-free graph $H$ such that $H$ is Ramsey for $\cT_{N^+_{B}(1)}$. Let $uv$ be an edge disjoint from $H$. We now consider two possible cases.\medskip

\noindent\textbf{Case 1}: $N_{D}(1)\subsetneq [q]\setminus \set{1}$.
Join $uv$ and every edge $h$ of $H$ by a copy of $T$, identifying the vertex $u$ with $a\in V(T)$, the vertex $v$ with $b \in V(T)$, and the edge $h$ with $cd \in E(T)$ (in an arbitrary fashion). We assume that, outside of the specified intersections, all copies of $T$ are vertex- and edge-disjoint. Call the resulting graph $G$. We claim that $G$ is an $X$-determiner with signal edge $uv$ for some non-empty $X\subsetneq [q]$.

Because $N_{D}(1)\subsetneq [q]\setminus \set{1}$, there exists $k\neq 1$ such that $k\not\in N_{D}(1)$; then $1\not\in N_D(k)$ and $1\in N_{D^c}(k)$. By \thref{lem:basic_conditions_new}, there exists $j\in N_D(k)$. Note that we trivially have $j\neq 1$ and $j\neq k$. Let $\phi$ be the following coloring of $G$. Set $\phi(uv) = j$. For every copy of $T$ in $G$, assign color $k$ to the copy of $bc$, and for any edge $h\in E(H)$, put $\phi(h)=1$. As $k\in N_D(j)$ and $1\in N_{D^c}(k)$, by~\thref{lem:SSc_properties_new} we can extend this coloring to the copies of $T$ so that each of these copies has a $\cT$-free coloring. We claim that this yields a $\cT$-free coloring of all of $G$. Suppose not; then there exists a color $\ell\in [q]$ and a copy $K$ of $K_{t_\ell}$ that is monochromatic in color $\ell$. If $K$ contains a vertex in a copy of $T$ that is not contained in $V(H)\cup \set{u,v}$, then this vertex does not have any neighbors outside of this copy of $T$. Then the clique $K$ must be fully contained in this copy of $T$, which is not possible by our choice of coloring. Thus, $V(K) \subseteq V(H)\cup \set{u,v}$. But the subgraph of $G$ induced by $V(H)\cup \set{u,v}$ consists of a single edge in color $j$, a star in color $k$, and the $K_{t_1}$-free graph $H$, which is monochromatic in color $1$, a contradiction since $t_\ell\geq 3$. 

Suppose $\phi$ is a $\cT$-free coloring of $G$ with $\phi(uv)=1$. Then every edge of $H$ is colored with a color in $N^+_{B}(1)$. As $H\to \cT_{N^+_{B}(1)}$, this coloring is not $\cT$-free. Therefore $G$ is an $X$-determiner for $\cT$, for some $X$ with $j\in X$ and $1\not\in X$, as required. By \thref{thm:all_det_are_safe}, it is furthermore a safe determiner. \medskip

\noindent\textbf{Case 2}: $N_{D}(1) = [q]\setminus \set{1}$, i.e., $N_{D^c}(1)=\set{1}$. This second case is very similar to the first one, using $T'$ instead of $T$. Recall that $H$ is a $K_{t_1}$-free graph such that $H\to \cT_{N^+_{B}(1)}$, and $uv$ is an edge disjoint from $H$. For each edge $h$ of $H$, take a copy of $T'$ and identify the copy of $pa$ in $T'$ with the edge $uv$ and the copy of $cd$ in $T'$ with the edge $h$ in an arbitrary fashion. Call the resulting graph $G$. As before, we will show that $G$ is an $X$-determiner for $\cT$ with signal edge $uv$ for some non-empty set $X$ with $1\not\in X$.

We begin by showing that $G\nto_q \cT$. For this, let $j\in N^-_B(1)$ and $k\in N_{D^c}(j)$. By~\thref{lem:SSc_properties_new}\ref{lem:SSc_colored_edges_new}, we have $1\not\in N^-_{B}(1)$, and therefore $j\neq 1$. As $N_{D^c}(1)=\set{1}$, we know that $j\in N_D(1)$, and hence $1\not\in N_{D^c}(j)$ and therefore $k\neq 1$. Then by~\thref{lem:SSc_properties_new}\ref{lem:ScSSc_any_color_new}, $T'$ has a $\cT$-free coloring in which the edges $pa$, $ab$ and $cd$ receive colors $k$, $j$, and $1$, respectively. We can then define a coloring $\phi$ of $G$ by setting $\phi(uv) = k$ and $\phi(h) = 1$ for every $h\in E(H)$, and extending this coloring to all copies of $T'$ so that each of these copies receives a $\cT$-free coloring. To verify that $\phi$ is a $\cT$-free coloring of $G$, notice that any clique containing a vertex from a copy of $T'$ that is not in $V(H)\cup \set{u,v}$ must be fully contained in that copy of $T'$. So it suffices to show that the subgraph of $G$ induced by $V(H)\cup \set{u,v}$ is $\cT$-free, which is easily seen to be true similarly to Case $1$. 
Note that, since $N_D(1) = [q]\setminus\set{1}$, we have $N_{D^c}(1) = \set{1}$ and thus $N^+_{B'}(1) = N^+_B(1)$ and $H\to \cT_{N^+_{B'}(1)}$. Therefore, as before, any coloring of $G$ that assigns the color $1$ to $uv$ cannot be $\cT$-free. 

Therefore $G$ is an $X$-determiner for $\cT$, for some $X$ with $k\in X$ and $1\not\in X$, as required. By \thref{thm:all_det_are_safe}, it is furthermore a safe determiner. 
\end{proof}\smallskip

We now put everything together and prove our main theorem.

\begin{proof}[Proof of \thref{thm:main_existence}]
We reason by strong induction on the number of colors $q$. For $q=2$, following~\cite{burr_graphs_1976,burr1977ramseyminimal} we know that for any $s,t\geq 3$ there exist signal senders for $(K_s,K_t)$ whose signal edges are arbitrarily far apart and that for any distinct integers $s,t\geq 3$, there exist red- and blue- determiners for $(K_s,K_t)$. By~\thref{thm:all_det_are_safe,thm:all_senders_are_safe} these gadgets are safe, and therefore any pair of cliques is distinguishable.

Let $q\geq 3$ and assume that any tuple of cliques with at most $q-1$ elements is distinguishable. Let $\cT=(K_{t_1},\ldots,K_{t_q})$ be a $q$-tuple of cliques with $t_1\geq \ldots\geq t_q\geq 3$. By \thref{lem:final_B_complete_new} there exists at least one safe set-determiner $\cT$, which together with \thref{thm:one_det_is_enough} concludes the induction proof. 
\end{proof}

\section{Conclusion}

We recall the following conjecture from~\cite{fox2016minimum} on the monotonicity of the minimum degree parameter.
\begin{conj}[\FGLPS{}~\cite{fox2016minimum}]
For all $q \geq 3$ and $t \geq 3$ we have $s_q(K_t) \geq s_q(K_{t-1})$.
\end{conj}
We believe that the gadgets for asymmetric tuples developed in this article could be useful to prove the conjecture. Indeed, it is a direct consequence of the following statement, now accessible using our new set-senders and set-determiners for asymmetric tuples of cliques.
\begin{conj}
For any $q\geq 2$, $t\geq 3$, and $0< k\leq q$, we have
\[s_q(\underbrace{K_t,\ldots,K_t}_{k \text{ times}},\underbrace{K_{t+1},\ldots,K_{t+1}}_{q-k \text{ times}}) \leq s_q(\underbrace{K_t,\ldots,K_t}_{k-1 \text{ times}},\underbrace{K_{t+1},\ldots,K_{t+1}}_{q-k+1 \text{ times}}).\]
\end{conj}

The use of signal senders and determiners in Ramsey theory extends beyond the setting of cliques. For instance, \Burr{}, \Nesetril{}, and \Rodl{}~\cite{burr1985useofsenders} showed the existence of appropriate gadgets for pairs $(H_1,H_2)$, where each of $H_1$ and $H_2$ is either 3-connected or isomorphic to $K_3$, and proved analogs of many of the results from~\cite{burr_graphs_1976,burr1977ramseyminimal} in this setting as applications. Versions of the aforementioned results of~\Rodl{} and~\Siggers{}~\cite{rodl_ramsey_2008} were also shown by \Siggers{} for odd cycles in~\cite{siggers_highly_2008,siggers_five_2012}, for pairs of non-bipartite 3-connected graphs~\cite{siggers_non-bipartite_2014}, and for certain pairs including cycles~\cite{siggers_non-bipartite_2014}.

In light of these results, the question of existence of set-senders and set-determiners for any tuple of $3$-connected graphs would be a natural continuation of the line of research pursued in this paper.
\begin{conj}
Any tuple $\cT=(H_1,\ldots,H_q)$, where $H_i$ is $3$-connected or isomorphic to $K_3$, is distinguishable.
\end{conj}

\section*{Acknowledgements} (SB) The research leading to these results was partially supported by the Deutsche Forschungsgemeinschaft (DFG, German Research 
Foundation) under Germany's Excellence Strategy – The Berlin Mathematics 
Research Center MATH+ (EXC-2046/1, project ID: 390685689, BMS Stipend) and Graduiertenkolleg ``Facets of Complexity'' (GRK 2434) and by EPSRC, grant no.\ EP/V048287/1. There are no additional data beyond that contained within the main manuscript.
(TL) The author was supported by the Commonwealth through an Australian Government Research Training Program Scholarship. 
We are grateful to Anita Liebenau for many helpful discussions. We also thank Andrzej \Rucinski{} for bringing the paper~\cite{graham2002ramsey} and the fact that \thref{thm:Ramsey_equiv} had been shown previously to our attention. We are also grateful to the anonymous referees for their detailed feedback. 

\bibliographystyle{abbrv}
\bibliography{Bibliography.bib}

\end{document}